\documentclass{article}

\newcommand{\dy}{\,\mathrm{d}\mathbf{y}}
\newcommand{\normal}{\hat{\mathbf{n}}}

\usepackage{graphicx}
\usepackage{placeins}
\usepackage{hyperref}
\usepackage{subcaption}
\usepackage{amsmath,amssymb,amsfonts}
\usepackage{multirow}
\usepackage{xcolor}
\usepackage{booktabs}
\usepackage{algorithm}
\usepackage[noEnd]{algpseudocodex}

\usepackage[style=numeric, sorting=none, maxnames=99]{biblatex}
\title{Robust Hierarchical Matrix Compression of Acoustic Volume and Boundary Integral Operators}
\author{Alberto Almuna-Morales$^{1,2}$ \and Danilo Aballay$^1$ \and Ignacio Labarca-Figueroa$^1$ \and Elwin van 't Wout$^1$\thanks{Contact: e.wout@uc.cl}}
\date{July 21, 2026}

\begin{document}

\maketitle

\begingroup
\footnotetext[1]{Institute for Mathematical and Computational Engineering, School of Engineering and Faculty of Mathematics, Pontificia Universidad Católica de Chile, Santiago, Chile.}
\footnotetext[2]{Department of Mechanical Engineering, University College London, London, United Kingdom.}
\endgroup

\begin{abstract}
Discretizing integral formulations of the Helmholtz equation yields dense linear systems. Hence, simulating acoustic models at larger scales or higher frequencies is typically constrained by memory capacity. Fast algorithms, such as hierarchical matrix compression, reduce the memory footprint substantially while controlling the approximation error in matrix-vector multiplications. However, the commonly used Adaptive Cross Approximation suffers from early-convergence problems, where the iterative construction of low-rank decompositions stops before reaching the targeted error tolerance. This failure arises when the error estimator does not capture significant components of the matrix structure under partial pivoting. This manuscript proposes a new diagonal convergence criterion, additional matrix elements for the pivoting strategy, an extended admissibility condition, and a sustained convergence check to improve the robustness of hierarchical matrix compression. These modifications improve compression reliability without increasing memory consumption. We tested our compression strategy on various discretized volume and boundary integral operators. The computational results show that our approach successfully compresses all benchmark matrices within predefined tolerances, thereby resolving the early-convergence issues encountered in standard algorithms. This robust matrix compression was achieved at the same memory footprint as alternative compression strategies. Furthermore, a complexity analysis shows log-linear memory scaling with mesh refinement at constant frequency. Finally, we successfully applied our robust matrix compression algorithm to a coupled system of volume and boundary integral operators that models transcranial ultrasound propagation. This confirms the feasibility of our robust algorithm to accelerate large-scale simulations with high-resolution meshes in a biomedical application.
\end{abstract}

\section{Introduction}
\label{sec:Introduction}

Performing realistic simulations in computational science and engineering requires fast, accurate, and robust numerical methods for solving partial differential equations. This manuscript presents numerical solutions of the Helmholtz equation for acoustics with coupled volume and boundary integral operators. As these operators involve large, dense matrices, fast and accurate simulations can be performed only with specialized compression techniques that enable memory-efficient matrix-vector multiplications. However, standard hierarchical low-rank approximations may lack robustness, as numerical errors may exceed the prescribed tolerances in certain scenarios.

This manuscript aims to mitigate the early-convergence problems of standard compression techniques to achieve robust calculations in acoustic integral equation methods. We will consider the Helmholtz equation, which is a linear partial differential equation for scalar harmonic waves that has been extensively used in physics and engineering~\cite{ihlenburg1998finite, lahaye2017modern}. Integral equations are well-established formulations of the Helmholtz equation for acoustic transmission models~\cite{nedelec2001acoustic, steinbach2008numerical}. They can be solved with numerical approaches such as boundary integral formulations via the Boundary Element Method (BEM)~\cite{sauter2010boundary, wout2021benchmarking}, and coupled volumetric-boundary frameworks using Volume-Surface Integral Equations (VSIE)~\cite{costabel2015spectrum, labarca2024volume}. These integral equation methods naturally handle the homogeneous, unbounded exterior domain, as the exterior acoustic field can be expressed in terms of a Green's function that satisfies the radiation condition~\cite{labarca2025coupled}. Hence, the computational domain is limited to material interfaces and bounded subdomains with localized heterogeneous material properties~\cite{aballay2026nested}. Integral equation methods are also efficient at simulating wave propagation at high frequencies~\cite{wout2022pmchwt, aubry2022benchmark}. However, these methods are memory-intensive, as the discretization matrix is dense.

To solve the systems of discrete integral operators, we employ the Generalized Minimal Residual Method (GMRES)~\cite{saad1986gmres, liesen2012krylov}. Given the dense nature of the discretization matrices, the most expensive routine in this iterative Krylov subspace method is the matrix-vector product. Fast calculation of the matrix-vector product with a significantly reduced memory footprint can be achieved with algorithms such as the Fast Multipole Method~\cite{greengard1987fast, darve2000fast, wang2021exafmm} and hierarchical matrix compression ($\mathcal{H}$-matrices)~\cite{borm2010efficient, bebendorf2008hierarchical}. The Fast Multipole Method decouples near- and far-field calculations via a tree-based hierarchy of analytic expansions. In contrast, hierarchical matrix compression leverages low-rank representations of matrix blocks originating from far-field interactions, while using dense blocks for the fewer near-field interactions.
This manuscript considers standard variants of $\mathcal{H}$-matrix compression frameworks, as their efficiency and performance are sufficient for many practical applications of acoustic scattering across low- and high-frequency regimes~\cite{wout2015fast, gelat2025evaluation}. They are also the building blocks for more advanced approaches such as $\mathcal{H}^2$-matrices~\cite{borm2010efficient} and directional compression~\cite{borm2017approximation, borm2024memory}.

Hierarchical matrix compression can achieve an $O(n \log n)$ storage complexity~\cite{betcke2017computationally}, where $n$ is the number of degrees of freedom (DOFs) in the linear system. This represents a substantial reduction from the $O(n^2)$ memory requirement of dense matrices. The guiding principle behind the compression strategy is that the Helmholtz Green's function permits accurate, separable approximations of the source and observer influences for distant interactions~\cite{engquist2018approximate}. Consequently, a hierarchical matrix structure is created with large blocks for far-field interactions and small blocks for near-field interactions between DOFs. The large matrix blocks are then approximated with a low-rank decomposition. The most commonly used decomposition algorithm is the Adaptive Cross Approximation (ACA)~\cite{grasedyck2005}. Its full-pivoting variant requires the full assembly of all blocks, yielding a prohibitive $\mathcal{O}(n^2)$ storage complexity. In contrast, partial-pivoting variants iteratively assemble matrix elements upon demand to construct low-rank approximations. It uses only a subset of each block's rows and columns, which is essential to maintain the desired log-linear complexity. For this reason, we will use Adaptive Cross Approximation with Partial Pivoting (ACAPP).

A major drawback of ACAPP is that partial-pivoting implementations can suffer from early-convergence problems. In those cases, the iterative creation of the low-rank decomposition terminates too early due to an inaccurate error estimate from incomplete matrix information. This phenomenon has been reported in the literature (e.g.,~\cite{Heldring2021, tetzner_2024, Laviada2009, Heldring2014, Heldring2015}). In~\cite{Heldring2021}, the authors propose a mean-based random sampling convergence criterion~(RSCC) to improve ACAPP convergence. In~\cite{tetzner_2024}, early-convergence issues were reported for the magnetic field integral equation for electromagnetics. As remedies, a combination of the standard criterion with the RSCC and geometry-based pivoting strategy was proposed. The authors in~\cite{Laviada2009} propose a multiple-restart approach to mitigate early-convergence issues. Finally, the ACA+ version leverages additional information in the pivoting strategy, thereby improving ACA convergence in domains with edges~\cite{grasedyck2005}.

While practical examples of early convergence have been well-documented for boundary integral operators, that is not the case for volume integral operators. Even though the basic ACA design is agnostic to the operator's dimension, the compression performance strongly depends on the specific integral operator. This manuscript shows that ACAPP also exhibits early convergence for volume integral operators. Moreover, mitigation strategies that worked for the BEM are insufficient to achieve robust compression of the VSIE.

The primary objective of this manuscript is to enhance the robustness of hierarchical matrix compression of the discretized volume and boundary integral operators for acoustics. Specifically, we will mitigate early-convergence issues encountered in standard ACAPP implementations. We propose four algorithmic extensions that are deterministic and preserve the log-linear storage complexity. First, we design a new diagonal convergence criterion based on information from the extended main diagonal of the block. Second, we augment the algorithm's pivoting strategy by incorporating the diagonal elements as candidate pivots. Third, we implement an admissibility condition that prevents near-singular interactions in the admissible blocks. Fourth, we develop a sustained convergence check for the error estimator. 

The numerical experiments demonstrate the advantages of our strategy. Importantly, the proposed methodology successfully resolves all encountered early-convergence problems while maintaining the efficiency of the original compression algorithms.

This manuscript is organized as follows. Section~\ref{sec:matrix_compression} introduces the methodology for standard $\mathcal{H}$-matrix compression. We present our algorithmic extensions to achieve robust ACAPP compression for acoustics in Section~\ref{sec:robust_compression}. The numerical results in Section~\ref{sec:results} demonstrate the robustness of our matrix compression algorithm for volume and boundary integral operators across several benchmarks, including canonical geometries and a CT-derived human skull model. 

\section{Matrix compression}
\label{sec:matrix_compression}

This section describes the methodology of hierarchical matrix compression of discretized volume and boundary integral operators, following standard literature (e.g., \cite{borm2010efficient, bebendorf2008hierarchical, hackbusch2004hierarchical, borm2003introduction}).

\subsection{Integral operators}

We consider time-harmonic acoustic wave propagation and linear material responses, using the time convention $e^{-\imath \omega t}$. The acoustic pressure $p(\mathbf{x})$ is then described by the heterogeneous Helmholtz equation
\begin{equation}
    \rho(\mathbf{x})\nabla\cdot\left( \frac{1}{\rho(\mathbf{x})} \nabla p(\mathbf{x})\right) + \left( \frac{\omega}{c(\mathbf{x})}\right)^2 p(\mathbf{x}) = 0, \quad\mathbf{x}\in\Omega;
    \label{eq:helmholtz}
\end{equation}
where $\Omega$ is a bounded, open domain. Here, $\omega$ denotes the angular frequency, $\rho(\mathbf{x})$ the density, and $c(\mathbf{x})$ the speed of sound. The boundary of $\Omega$ is denoted by $\Gamma$ and equipped with a unit outward-pointing normal vector $\normal$. The unbounded region exterior to the surface is considered homogeneous with constant density and speed of sound, denoted by $\rho_0$ and $c_0$, respectively. The Green’s function associated with the three-dimensional Helmholtz equation for constant wavenumber $\kappa = \omega/c_0$ is given by
\begin{equation}
    G_{\kappa}(\mathbf{x}, \mathbf{y})=\frac{1}{4 \pi} \frac{e^{\imath \kappa|\mathbf{x}-\mathbf{y}|}}{|\mathbf{x}-\mathbf{y}|}\quad \text{for} \quad \mathbf{x} \neq \mathbf{y},
    \label{eq:green_function_helmholtz}
\end{equation} 
where $\imath$ denotes the imaginary unit~\cite{duffy2015green}.

The BEM for piecewise-constant material parameters and the VSIE for local heterogeneous materials are different numerical methods, but share many characteristics. For example, both share the same type of integral operators. Specifically, we will consider the single- and (adjoint) double-layer integral operators. Many boundary integral formulations can be created with these three operators~\cite{wout2021benchmarking}. Although the proposed matrix compression algorithms can also be applied to the hypersingular integral operator, we will not consider it here for brevity. Precisely, the single-layer volume integral operator (SL-VIO) is given by
\begin{equation}
S\left[\phi\right](\mathbf{x}) = \iiint_{\Omega} G_{\kappa}(\mathbf{x},\mathbf{y}) \frac{\rho_0}{\rho(\mathbf{y})} \left(\left(\frac{\omega}{c(\mathbf{y})}\right)^2 - \left(\frac{\omega}{c_0}\right)^2\right) \phi(\mathbf{y}) \dy, \quad\mathbf{x}\in\Omega.
\label{eq:sl-vio}
\end{equation}
The double-layer boundary integral operator (DL-BIO) can be written as 
\begin{equation}
K \left[\phi\right](\mathbf{x})=\iint_{\Gamma} \frac{\partial}{\partial \hat{\mathbf{n}}(\mathbf{y})} G_{\kappa}(\mathbf{x}, \mathbf{y}) \phi(\mathbf{y}) \dy, \quad\mathbf{x}\in\Gamma.
\label{eq:dl-bio}
\end{equation}
Finally, the adjoint double-layer volume integral operator (ADL-VIO) reads
\begin{equation}
T \left[\phi\right](\mathbf{x}) = \nabla_{\mathbf{x}} \cdot \iiint_{\Omega} G_{\kappa}(\mathbf{x},\mathbf{y}) \, \nabla_{\mathbf{y}} \left( \frac{\rho_0}{\rho(\mathbf{y})}\right) \phi(\mathbf{y}) \dy, \quad\mathbf{x}\in\Omega.
\label{eq:adl-vio}
\end{equation}
See, e.g., \cite{sauter2010boundary, wout2021benchmarking} and \cite{costabel2015spectrum, labarca2024volume, labarca2025coupled, aballay2026nested} for further details about boundary and volume integral operators, respectively.

\subsection{Numerical discretization}

\begin{figure}[htbp]
    \begin{subfigure}[b]{0.49\linewidth}
        \centering
        \includegraphics[width=0.65\linewidth]{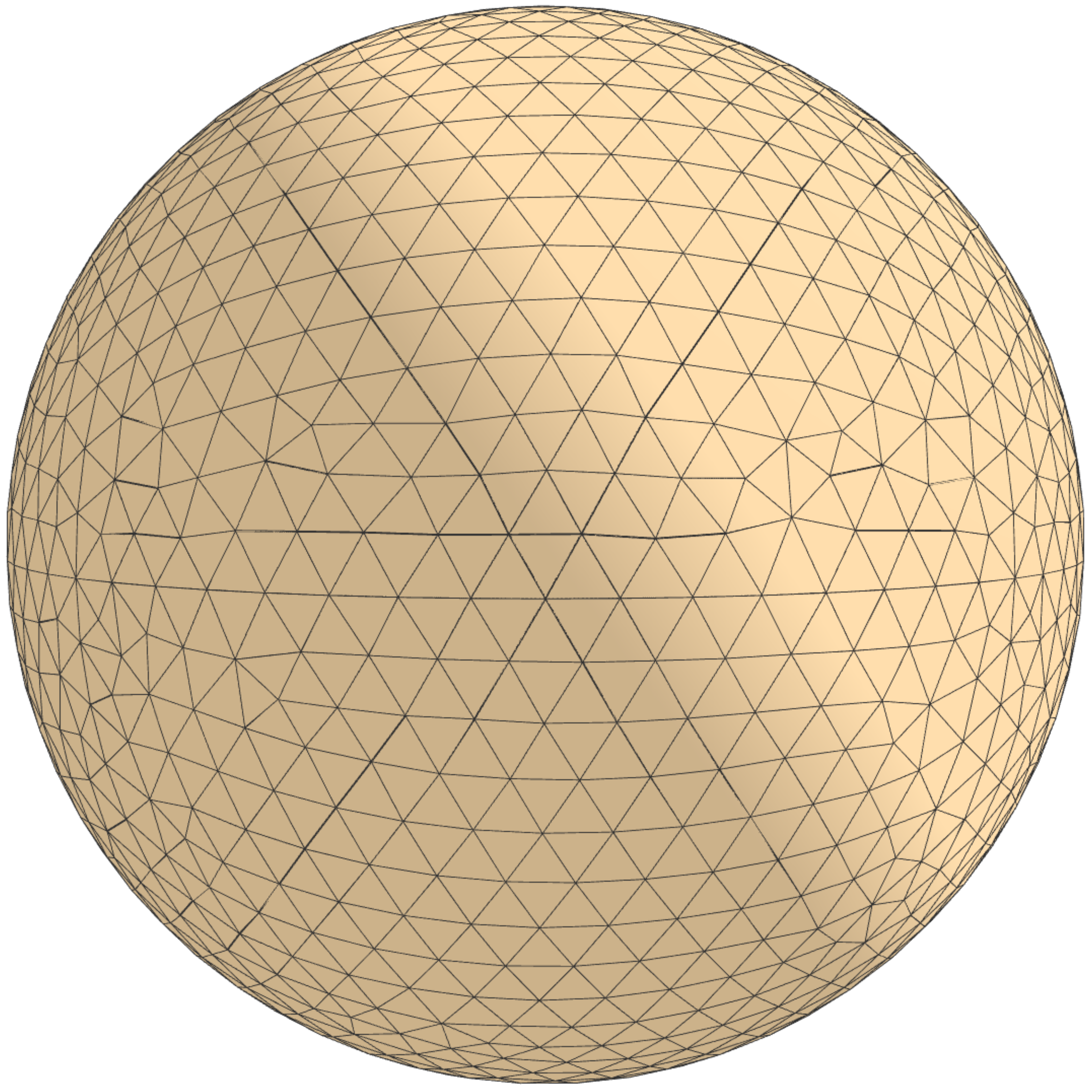}
        \caption{Triangular surface mesh.}
    \end{subfigure}
    \hfill
    \begin{subfigure}[b]{0.49\linewidth}
        \centering
        \includegraphics[width=0.65\linewidth]{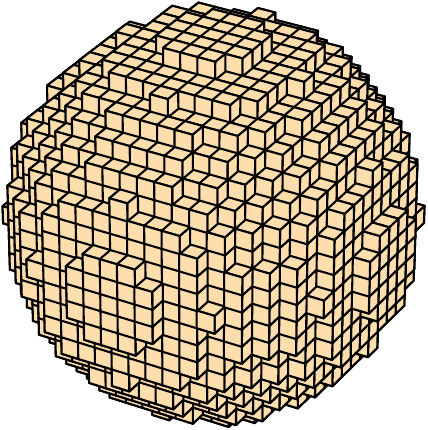}
        \caption{Voxel mesh.}
    \end{subfigure}
    \caption{Example meshes on a sphere.}
    \label{fig:sphere_meshes}
\end{figure}

We will consider two different numerical discretization approaches. First, we consider a Galerkin discretization of boundary integral operators with piecewise linear (P1) test and basis functions on a triangular surface mesh. Second, we discretize volume integral operators using collocation and piecewise constant (P0) basis functions on a voxel mesh. This P0 collocation approach will also be applied to the boundary integral operators on the surface of the voxel mesh. The voxel mesh consists of rectangular cuboid elements that may have different resolutions in each axis. Hence, the DOFs for the P1 discretization are associated with the vertices of triangles, and the DOFs for the P0 discretization are associated with the centers of voxels or their faces at the surface. See Figure~\ref{fig:sphere_meshes} for example meshes.

\subsection{Octree partitioning}

We use an octree partitioning of the geometry to separate far from near interactions. We start with a bounding box around the geometry, which is then recursively divided into halves along each axis. Hence, each box is partitioned into eight smaller boxes, and the process continues until a predefined maximum tree depth or minimum box size is reached. Additionally, empty boxes with no DOFs are removed from the octree data structure.

\begin{figure}[htbp]
    \begin{subfigure}[b]{0.49\linewidth}
        \centering
        \includegraphics[width=\linewidth]{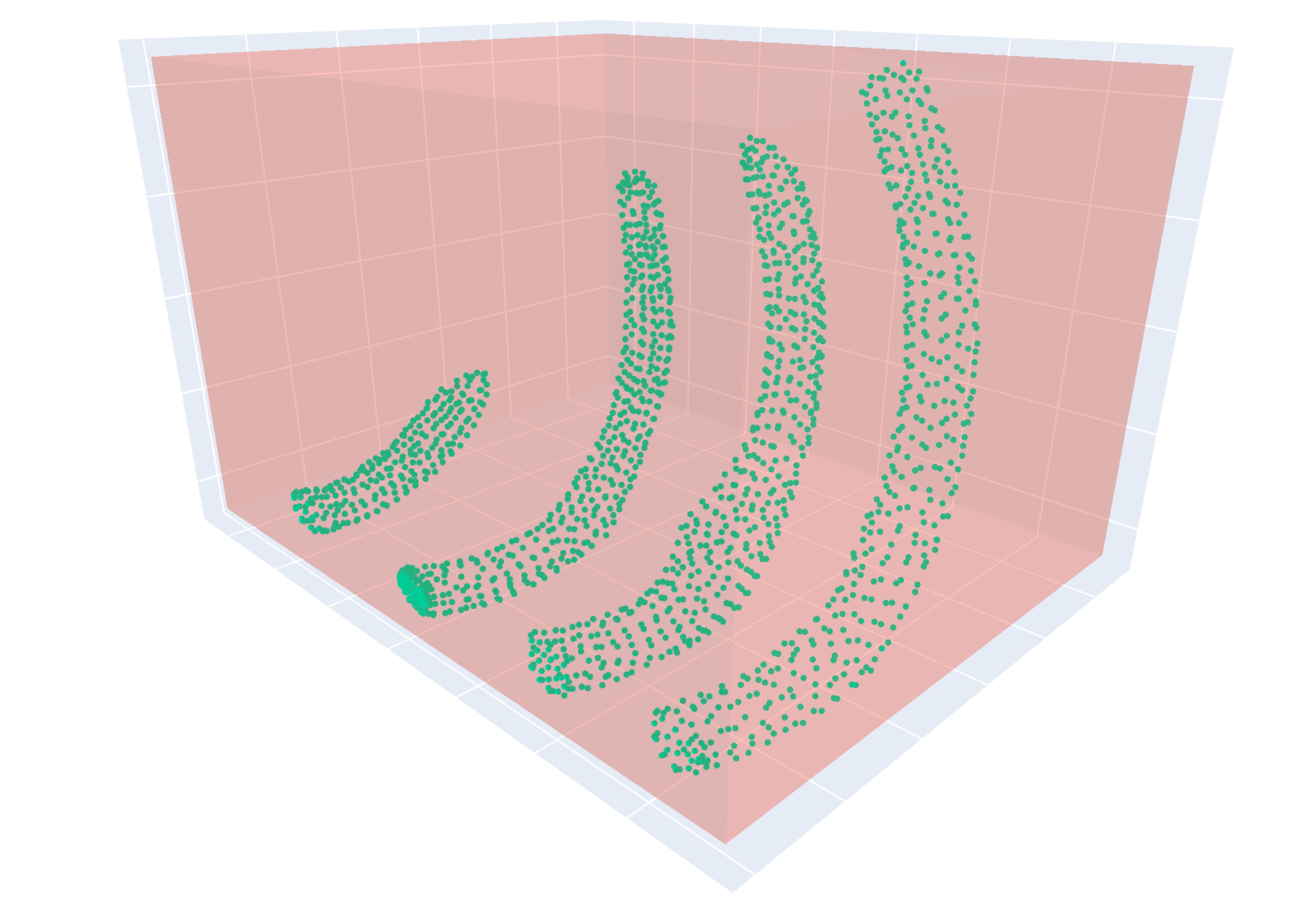}
    \end{subfigure}
    \hfill
    \begin{subfigure}[b]{0.49\linewidth}
        \centering
        \includegraphics[width=\linewidth]{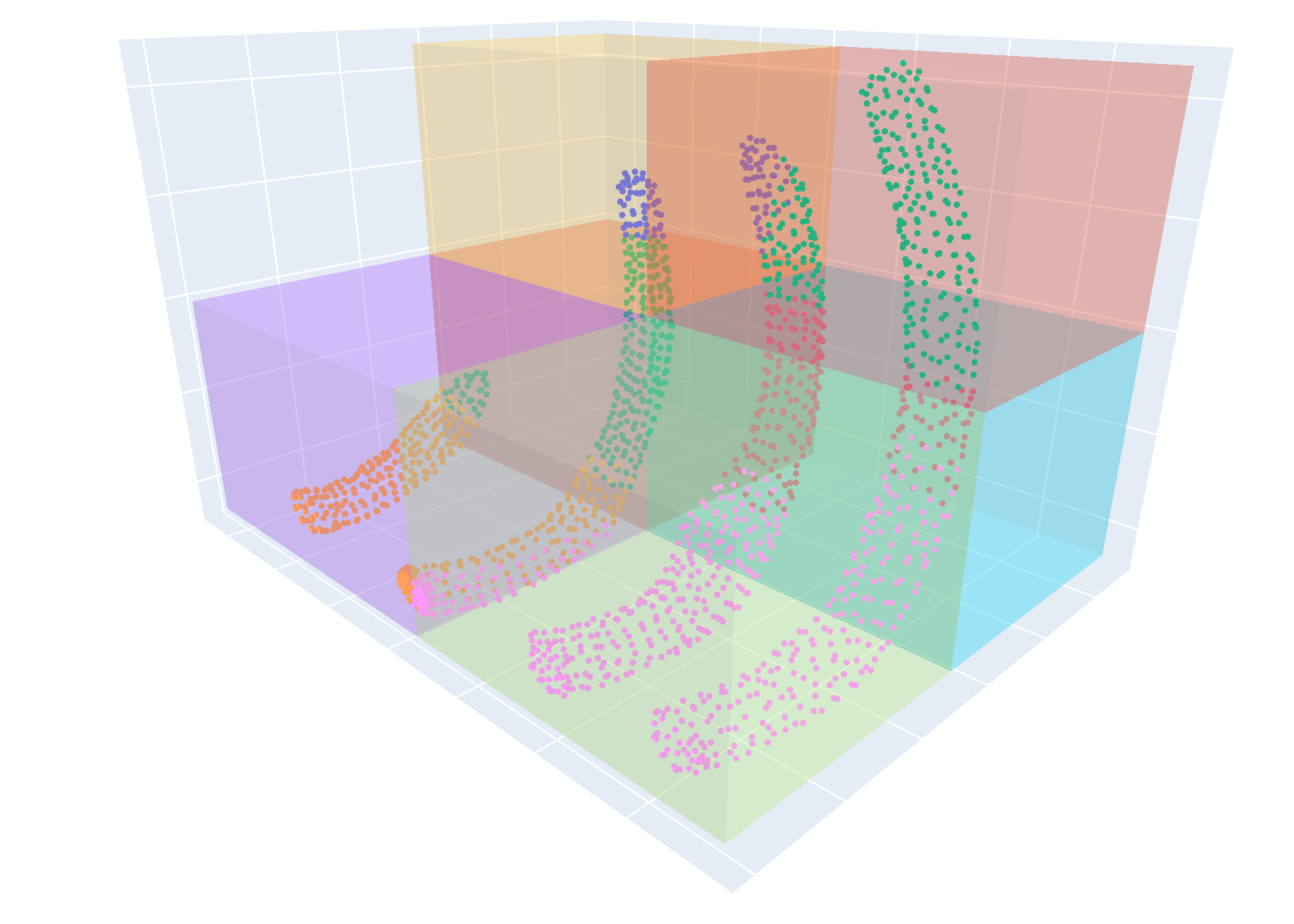}
    \end{subfigure}
    \begin{subfigure}[b]{0.49\linewidth}
        \centering
        \includegraphics[width=\linewidth]{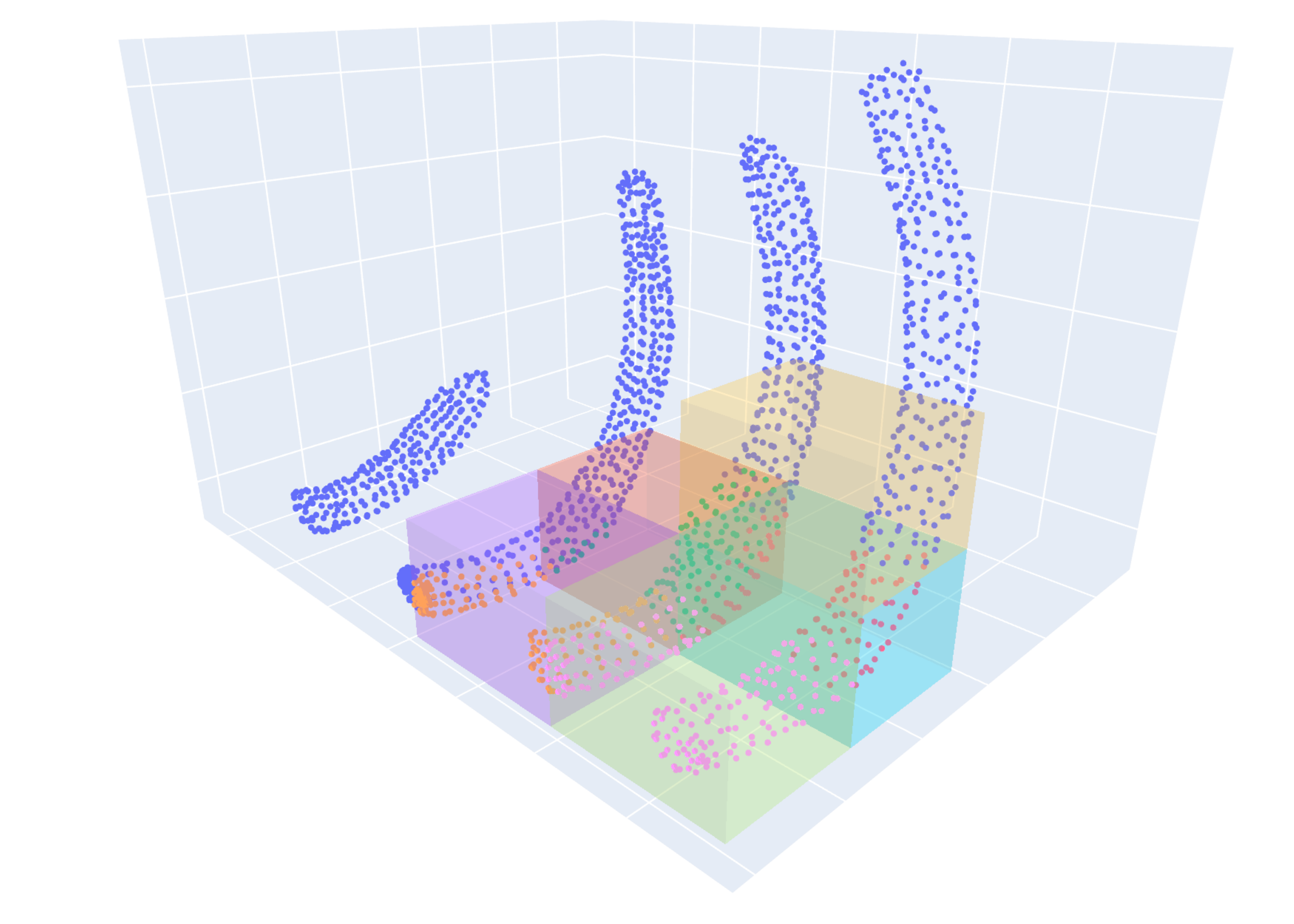}
    \end{subfigure}
    \hfill\begin{subfigure}[b]{0.49\linewidth}
        \centering
        \includegraphics[width=\linewidth]{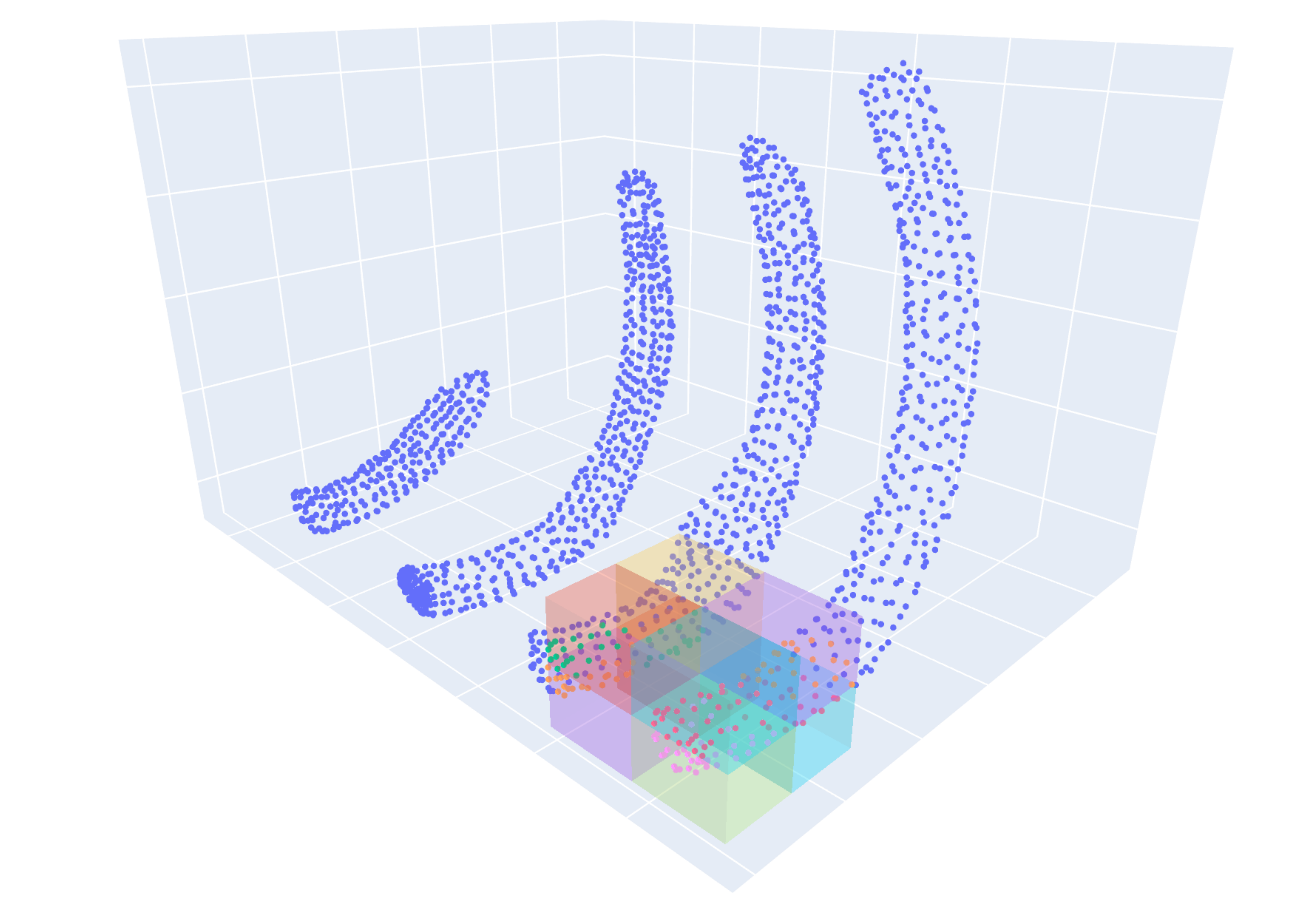}
    \end{subfigure}
    \caption{An example of the first four levels of a branch in an octree partitioning. The dots represent the vertices of a surface mesh of four ribs. Each colored volume represents a nonempty box in the octree.}
    \label{fig:branch_octree}
\end{figure}

The maximum tree depth and minimum box size are parameters that can be set in our implementation. The tree depth is the number of levels in the octree partitioning, while the box size indicates the number of DOFs within an octree box. These conditions avoid small matrix blocks that are unlikely to be efficiently represented by low-rank approximations, and their parameters can be set individually for each simulation to obtain small performance gains. Figure~\ref{fig:branch_octree} shows an example of an octree branch at four different levels.

\subsection{Admissibility}

\begin{figure}[htbp]
    \centering
    \includegraphics[width=.49\linewidth]{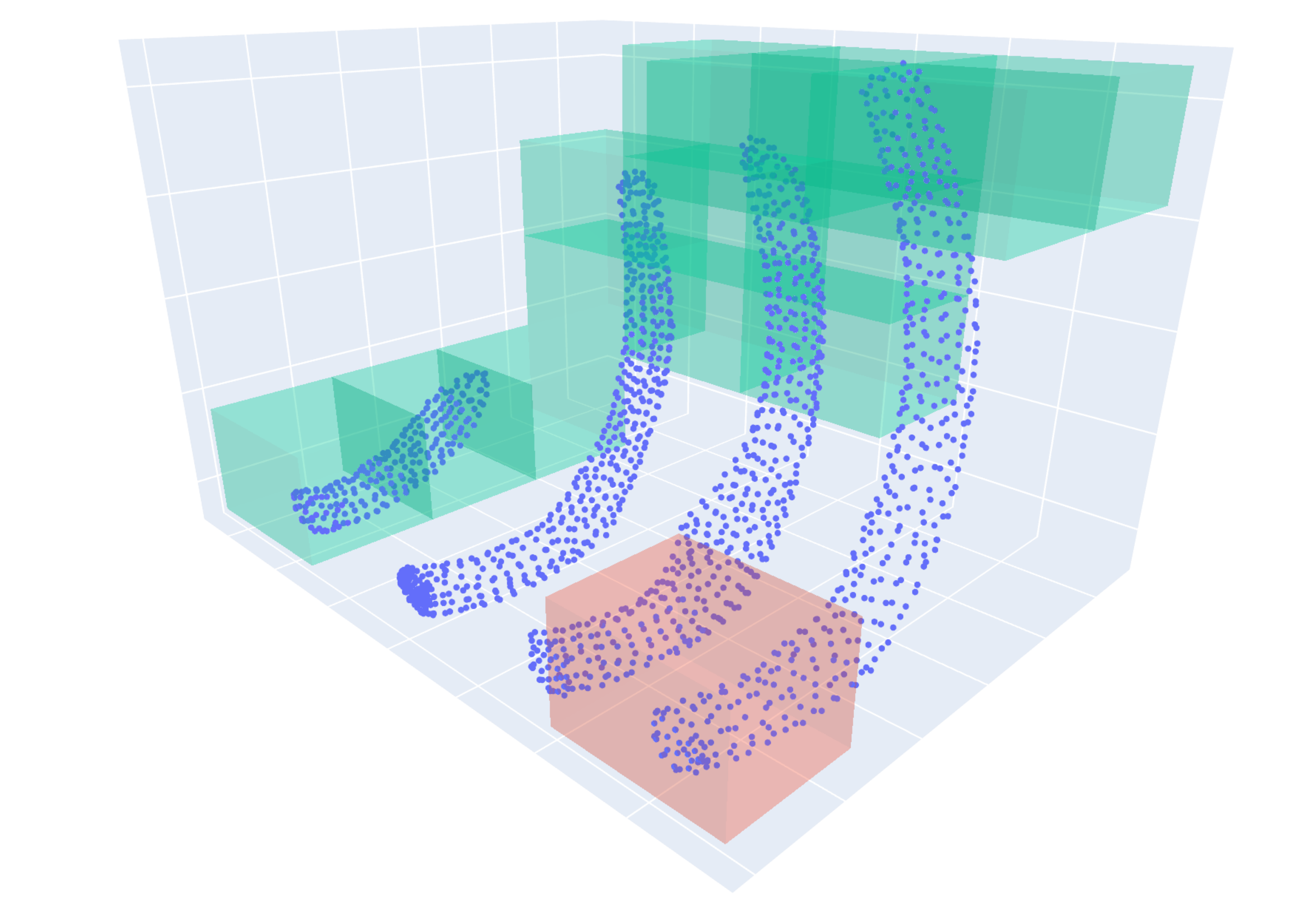}
    \includegraphics[width=.49\linewidth]{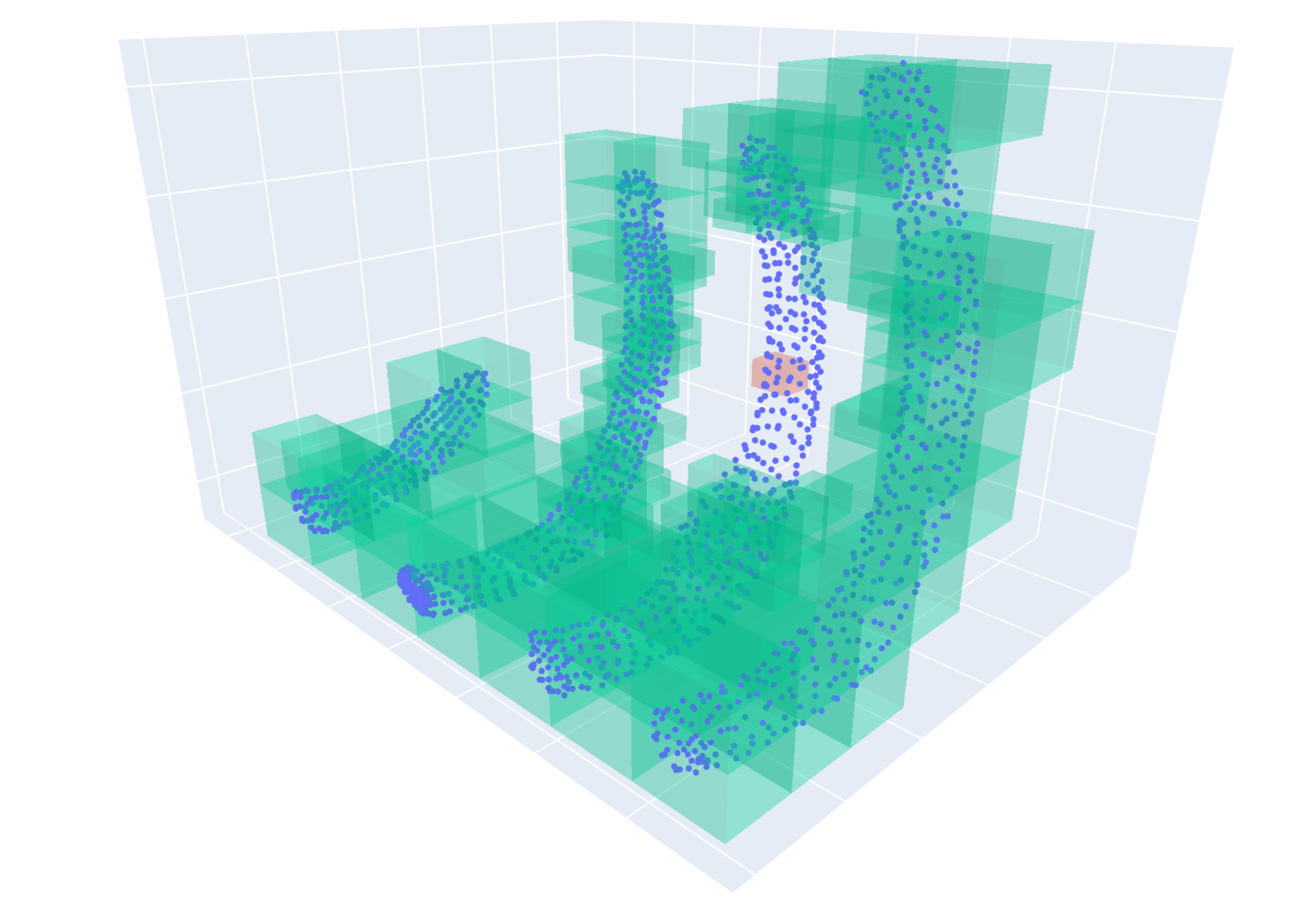}
    \caption{Examples of admissibility between boxes in an octree partitioning. The green boxes represent all admissible interactions with respect to the red one.}
    \label{fig:admissibility}
\end{figure}

The working principle of hierarchical matrix compression for integral equations is that geometrically separated regions yield matrix blocks that can efficiently be represented by low-rank approximations. To determine if a pair of boxes is considered a far or near interaction, we use the admissibility condition
\begin{equation}
    \operatorname{Adm}(Box_i, Box_j): \eta\operatorname{dist}(Box_i, Box_j) \geqslant \max(\operatorname{diam}(Box_i), \operatorname{diam}(Box_j)).
    \label{eq:admissibility}
\end{equation}
Here, the distance between two boxes is calculated as the shortest distance between all possible pairs of their corners, while the diameter of a box is the distance between its opposite corners. Furthermore, $\eta \in \mathbb{R}$ is a separation parameter. See Figure~\ref{fig:admissibility} for an illustration of the admissibility condition~\eqref{eq:admissibility} applied to the octree partitioning of a mesh.

Condition~\eqref{eq:admissibility} is one of the most common choices to assess admissibility, and is also known as a \emph{strong} admissibility condition~\cite{hackbusch2004hierarchical}. Other variants include taking the minimum rather than the maximum, and treating all separated boxes as admissible~\cite{hackbusch2004hierarchical, borm2010efficient}. In our computational experience, using the strong admissibility condition with $\eta = 1$ yielded robust and efficient compression of the Helmholtz integral operators.

\subsection{Low-rank approximation}

\begin{figure}[htbp]
    \centering
    \includegraphics[width=.7\linewidth]{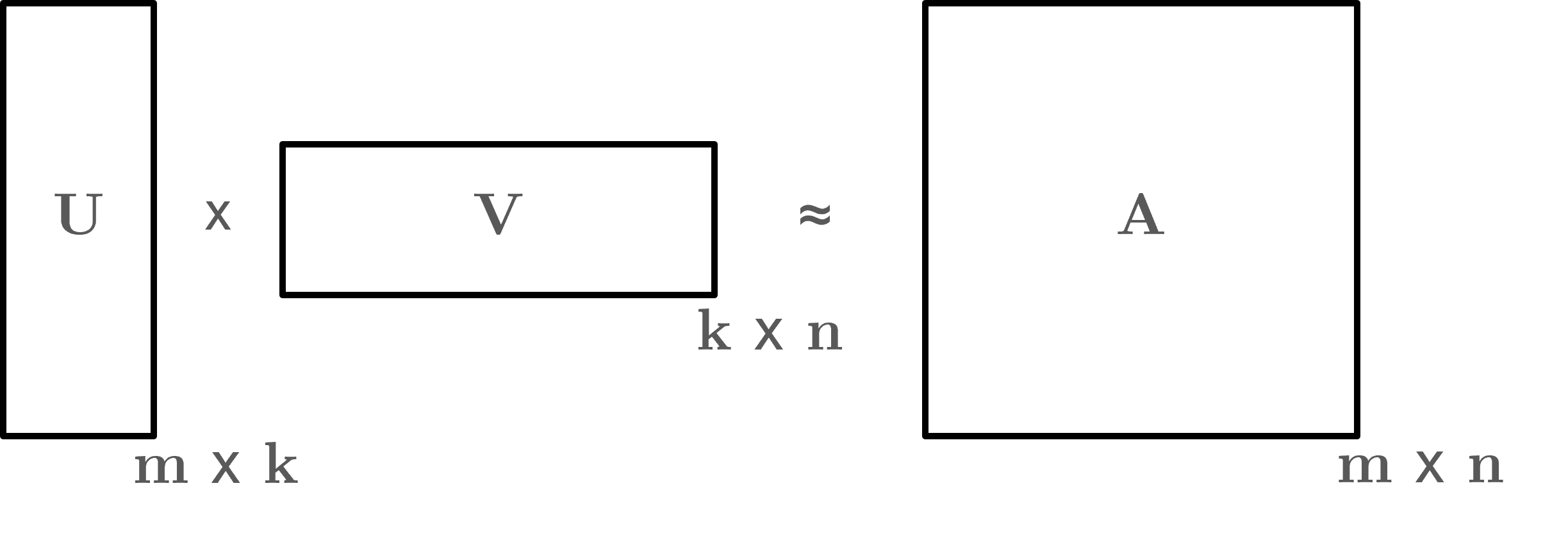}
    \caption{Low-rank approximation of rank $k$.}
    \label{fig:low-rank-approximation}
\end{figure}

The admissibility between two boxes in the octree determines how the integral operators acting on them will be assembled. Specifically, non-admissible interactions are always discretized into dense matrix blocks, while the admissible blocks are candidates for compression. Precisely, we search for a low-rank approximation $UV \approx A$ for matrix block $A \in \mathbb{C}^{m \times n}$ and low-rank matrices $U \in \mathbb{C}^{m \times k}$ and $V \in \mathbb{C}^{k \times n}$; see Figure~\ref{fig:low-rank-approximation}.

The goal of the matrix compression is to obtain an accurate decomposition for low ranks, i.e., $k \ll \min\{m,n\}$. We choose one of the most popular algorithms to achieve this low-rank approximation, namely the ACAPP~\cite{grasedyck2005}. Importantly, this heuristic algorithm does not require full access to the matrix block $A$ but instead evaluates specific rows and columns to iteratively build the low-rank representation. The ACAPP algorithm terminates when a compression error estimate falls below a predefined threshold, denoted by $\varepsilon$.

Algorithm~\ref{alg:acapp} shows an implementation of the default ACAPP. The notation for indexing vectors and matrices is based on the Python programming language. Notice that even though the matrix $A$ is an input parameter of the algorithm, it does not need to be explicitly constructed. In fact, each call to $A[i_\star,:]$ or $A[:,j_\star]$ is responsible for assembling a single row or column, respectively.

Our implementation of the ACAPP follows the standard algorithm with one extension. When the rank~$k$ becomes larger than half the original dimension, the low-rank approximation needs more memory to store $U$ and $V$ than the entire $A$. We detect this exception and switch to dense mode instead.

\begin{algorithm}[!ht]
\caption{Adaptive Cross Approximation with Partial Pivoting}\label{alg:acapp}
\begin{algorithmic}[1]

\Require $A\in\mathbb{C}^{m\times n}, \varepsilon\geqslant0$
\State $\left(m,n\right)\gets A.\texttt{shape}$ \Comment{Number of rows and columns of the matrix}
\State $i_\star\gets 0$ \Comment{Arbitrary selection of first row pivot}
\State $k\gets 0$ \Comment{Iteration number and compression rank}
\State $U^T \gets \texttt{[]}$ \Comment{Column vectors empty list}
\State $V \gets \texttt{[]}$ \Comment{Row vectors empty list}
\State $I\gets\{1,\dots, m\}$ \Comment{Available row pivots set}
\State $J\gets\{1,\dots, n\}$ \Comment{Available column pivots set}

\While{$\texttt{True}$}
    \State $\mathbf{v} \gets A[i_\star,:]$ \Comment{Take $i_\star$-th row of $A$}
    \State $I.\texttt{remove}(i_\star)$
    
    \State $\mathbf{v}\gets \mathbf{v} - \sum_{\ell=0}^{k-1}U^T[\ell,i_\star]\cdot V[\ell,:]$
    
    \State $j_\star\gets \texttt{argmax}_{j\in J}(|\mathbf{v}[j]|)$\Comment{Get column pivot}\label{alg_line:ACAPP_column_pivot}
    \State $\delta\gets \mathbf{v}[j_\star]$

    \If{$\delta == 0$}
        \If{$\texttt{length}(I)==0$} \Comment{All rows already checked?}
        \State \texttt{break} \Comment{Exit while}
        \EndIf
        \State $i_\star \gets \texttt{Random}(i\in I)$ \Comment{New random row pivot} \label{alg_line:random_step}
        \State $\texttt{continue}$ \Comment{Go to first line of while}
    \EndIf
    \State $\mathbf{v}\gets \mathbf{v} / \delta$
    \State $\mathbf{u} \gets A[:,j_\star]$ \Comment{Take $j_\star$-th column of $A$}
    \State $J.\texttt{remove}(j_\star)$
    
    \State $\mathbf{u}\gets \mathbf{u} - \sum_{\ell=0}^{k-1}V[\ell,j_\star]\cdot U^T[\ell,:]$
    
    \State $k\gets k+1$ \label{alg_line:aca_flag1}
    \State $U^T.\texttt{append}(\mathbf{u})$
    \State $V.\texttt{append}(\mathbf{v})$

    \State $\texttt{error}\gets \dots$ \Comment{Estimate error (Section~\ref{sec:error-estimation})} \label{alg_line:error_estimation}
    
    \If{$\texttt{error}\leqslant\varepsilon$} \Comment{Reached compression error?}
        \State \texttt{break} \Comment{Exit while}  \label{alg_line:aca_flag2}
    \EndIf 

    \State $i_\star\gets \texttt{argmax}_{i\in I}(|\mathbf{u}[i]|)$\Comment{Get row pivot}\label{alg_line:ACAPP_row_pivot}
\EndWhile
\State \Return $U, V$

\end{algorithmic}
\end{algorithm}

\subsection{Error estimation}\label{sec:error-estimation}

The iterative construction of the low-rank approximation $A \approx UV$ terminates when it reaches a predefined threshold $\varepsilon$. When the matrix block $A$ is available, the approximation error $\lVert A - UV \rVert$ can be calculated explicitly in an appropriate norm. However, the purpose of matrix compression is to reduce memory consumption. Hence, $A$ must not be assembled in its entirety, and the algorithm can have partial access to $A$ only to construct $U$ and $V$. This also means that the error must be estimated from partial information on the matrix blocks.

The error calculation on Line~\ref{alg_line:error_estimation} of Algorithm~\ref{alg:acapp} is left open to allow for custom choices. Since this is an essential design choice, let us explain one of the most common error estimators in detail. The ACAPP algorithm returns a $k$ low-rank approximation $R_k \approx A$ in the form of
\begin{equation}
    R_k=\sum_{\ell=1}^k \mathbf{u}_{\ell}\left(\mathbf{v}_{\ell}\right)^T=UV,
\end{equation}
where $\mathbf{u}_{\ell}$ and $\mathbf{v}_{\ell}$ are the column and row vectors calculated in iterations $\ell = 1, 2, \dots, k$. Ideally, the algorithm terminates when the relative approximation error satisfies
\begin{equation}
    \|A-R_k\| \leqslant \varepsilon\|A\|.
    \label{eq:ideal_inequality}
\end{equation}
However, the matrix $A$ cannot be accessed in its entirety. Therefore, the right-hand side will be replaced by $R_k$, obtaining
\begin{equation}
    \|A-R_k\| \leqslant \varepsilon\|R_k\|.
    \label{eq:approximate_inequality}
\end{equation}
On the left-hand side, we use the approximation
\begin{equation}
    \left\|A-R_k\right\| \lesssim\left\|A-R_{k-1}\right\| \approx\left\|R_k-R_{k-1}\right\|=\|\mathbf{u}_k\left(\mathbf{v}_k\right)^T\|,
    \label{eq:error_approximation}
\end{equation}
where $\|\mathbf{u}_k\left(\mathbf{v}_k\right)^T\| = \|\mathbf{u}_k\|_2\|\mathbf{v}_k\|_2$ for both the Frobenius norm and the 2-norm, as the matrix formed by the outer product is of rank 1. Combining the approximations yields the error estimation
\begin{equation}
    \frac{\|\mathbf{u}_k\|_2\|\mathbf{v}_k\|_2}{\|R_k\|} \leqslant \varepsilon.
    \label{eq:final_inequality}
\end{equation}
We evaluate the denominator with the Frobenius norm as
\begin{equation}
\left\|R_k\right\|_{\mathrm{F}}^2=\left\|R_{k-1}\right\|_{\mathrm{F}}^2+2 \operatorname{Re}\left(\sum_{{\ell}=1}^{k-1}\left(\mathbf{u}_{\ell}^* \mathbf{u}_k\right)\left(\mathbf{v}_{\ell}^* \mathbf{v}_k\right)\right)+\left\|\mathbf{u}_k\right\|_2^2\left\|\mathbf{v}_k\right\|_2^2,
\label{eq:fronorm_approximation}
\end{equation}
which can be updated during each iteration. Here, the star symbol $^*$ in $\mathbf{u}_{\ell}^*$ denotes the complex conjugate transpose of $\mathbf{u}_{\ell}$.

While condition~\eqref{eq:final_inequality} is a standard choice~\cite{tetzner_2024}, an alternative error estimate can be
\begin{equation}
    \frac{\|\mathbf{u}_k\|_2\|\mathbf{v}_k\|_2}{\|\mathbf{u}_1\|_2\|\mathbf{v}_1\|_2} \leqslant \varepsilon,
    \label{eq:euclidean_relative_error}
\end{equation}
which uses Euclidean norms~\cite{grasedyck2005}.

In our implementation, both conditions~\eqref{eq:final_inequality} and~\eqref{eq:euclidean_relative_error} must be satisfied for the ACAPP compression to terminate. Computational experience demonstrated increased robustness and minimal memory overhead when both termination conditions were required.

\subsection{Pivoting strategy}

A fundamental aspect of ACAPP algorithms is their pivoting strategy, which corresponds to lines \ref{alg_line:ACAPP_column_pivot} and \ref{alg_line:ACAPP_row_pivot} in Algorithm~\ref{alg:acapp}. The pivot in ACAPP is the matrix element that determines which column or row is calculated subsequently. Normally, the pivot for the next column or row is selected as the largest element in absolute value of the most recently updated row or column, respectively.

It has been observed that the ACAPP may fail for the double-layer operator in domains with edges~\cite{grasedyck2005}. The ACA+ uses an enhanced pivoting strategy that solves this specific setting~\cite{grasedyck2005}. As it also uses partial pivoting, we denote it by ACAPP+. See Algorithm~\ref{alg:aca_plus} for the implementation details. 

The main idea of ACAPP+ is to create a broader collection from which pivots can be selected. Precisely, two additional vectors, one matrix row and one column, are calculated but not directly incorporated in the low-rank approximation. The algorithm alternately selects pivots from these observer vectors and from the latest updated vector in the low-rank approximation. This strategy adds more diversity to the pivot search, exploring regions of the matrix block that the original ACAPP pivoting strategy would not prioritize.

\begin{algorithm}[!ht]
\caption{ACAPP+}\label{alg:aca_plus}
\begin{algorithmic}[1]

\Require $A\in\mathbb{C}^{m\times n}, \varepsilon\geqslant0$
\State $\left(m,n\right)\gets A.\texttt{shape}$ 
\State $\left(I, J\right)\gets\left(\{1,\dots, m\}, \{1,\dots, n\}\right)$
\State $j_{\text{ref}}\gets 0$ \Comment{Arbitrary selection of reference column pivot}
\State $\mathbf{u}^\text{ref}\gets A[:,j_{\text{ref}}]$ \Comment{Get reference column}
\State $i_{\text{ref}}\gets \texttt{argmin}_{i\in I}(|\mathbf{u}^\text{ref}[i]|)$\Comment{Get reference row pivot}
\State $\mathbf{v}^\text{ref}\gets A[i_{\text{ref}},:]$ \Comment{Get reference row}
\State $k\gets 0$
\State $\left(U^T, V\right) \gets \left(\texttt{[]}, \texttt{[]}\right)$
\While{$\texttt{True}$}
    \State $(i_\star, j_\star) \gets \left( \texttt{argmax}_{i\in I}(|\mathbf{u}^\text{ref}[i]|), \texttt{argmax}_{j\in J}(|\mathbf{v}^\text{ref}[j]|) \right)$ \label{alg_line:aca+_pivot_selection_start}
    \If{$|\mathbf{u}^\text{ref}[i_\star]|\geqslant|\mathbf{v}^\text{ref}[j_\star]|$}
        \State $\mathbf{v} \gets A[i_\star,:] - \sum_{\ell=0}^{k-1}U^T[\ell,i_\star]\cdot V[\ell,:]$
        \State $I.\texttt{remove}(i_\star)$
        \State $j_\star \gets \texttt{argmax}_{j\in J}(|\mathbf{v}[j]|)$
        \State $\delta\gets \mathbf{v}[j_\star]$
        \If{$\delta == 0$}
            \If{$\texttt{length}(I)==0$} 
            \State \texttt{break} 
            \EndIf
            \State $\texttt{continue}$ 
        \EndIf
        \State $\mathbf{u} \gets \left(A[:,j_\star] - \sum_{\ell=0}^{k-1}V[\ell,j_\star]\cdot U^T[\ell,:]  \right) \big/ \delta$
        \State $J.\texttt{remove}(j_\star)$
    \Else
        \State \# Analogous process to previous \textbf{if} statement (start with $\mathbf{u}$ then $\mathbf{v}$).
    \EndIf
    \State \# Lines \ref{alg_line:aca_flag1}$\rightarrow$\ref{alg_line:aca_flag2} of Algorithm~\ref{alg:acapp}.
    
    \State $\mathbf{u}^\text{ref}\gets \mathbf{u}^\text{ref} - \mathbf{v}[j_\text{ref}] \cdot \mathbf{u}$ \Comment{Update reference column} \label{alg_line:update_reference_start}
    \State $\mathbf{v}^\text{ref}\gets \mathbf{v}^\text{ref} - \mathbf{u}[i_\text{ref}] \cdot \mathbf{v}$ \Comment{Update reference row}
    \If{$i_\text{ref}==i_\star$ and $j_\text{ref}==j_\star$} \Comment{New reference column \& row}
        \State $j_\text{ref} \gets \texttt{Random}(j\in J)$\label{alg_line:aca+_random_reference} 
        \State $\mathbf{u}^\text{ref} \gets A[:,j_\text{ref}] - \sum_{\ell=0}^{k-1}V[\ell,j_\text{ref}]\cdot U^T[\ell,:]$
        \State $i_{\text{ref}}\gets \texttt{argmin}_{i\in I}(|\mathbf{u}^\text{ref}[i]|)$
        \State $\mathbf{v}^\text{ref} \gets A[i_\text{ref},:] - \sum_{\ell=0}^{k-1}U^T[\ell,i_\text{ref}]\cdot V[\ell,:]$
    \ElsIf{$i_\text{ref}==i_\star$} \Comment{New reference row}
        \State $i_{\text{ref}}\gets \texttt{argmin}_{i\in I}(|\mathbf{u}^\text{ref}[i]|)$
        \State $\mathbf{v}^\text{ref} \gets A[i_\text{ref},:] - \sum_{\ell=0}^{k-1}U^T[\ell,i_\text{ref}]\cdot V[\ell,:]$
    \ElsIf{$j_\text{ref}==j_\star$} \Comment{New reference column}
        \State $j_{\text{ref}}\gets \texttt{argmin}_{j\in J}(|\mathbf{v}^\text{ref}[j]|)$
        \State $\mathbf{u}^\text{ref} \gets A[:,j_\text{ref}] - \sum_{\ell=0}^{k-1}V[\ell,j_\text{ref}]\cdot U^T[\ell,:]$ \label{alg_line:update_reference_end}
    \EndIf
\EndWhile
\State \Return $U, V$
\end{algorithmic}
\end{algorithm}

\FloatBarrier

\section{Robust compression}
\label{sec:robust_compression}

While the standard ACAPP and ACAPP+ algorithms effectively provide low-rank approximations in most cases, they still suffer from early-convergence problems in various important situations. Here, we propose four improvements to achieve robust matrix compression.

\subsection{Early-convergence problem}
\label{sec:early-convergence}

The creation of low-rank approximations is an iterative process that continues adding vectors to the decomposition until a termination criterion is reached. However, as the ACAPP uses partial information about the matrix block, the termination criterion is not exact but instead relies on error estimates. Hence, the algorithm cannot guarantee that the actual compression error remains below the prescribed tolerance. In practice, the error estimation is sufficiently accurate to obtain a compressed matrix that satisfies the approximation threshold. However, various cases where the error estimation fails to deliver an effective low-rank approximation have been reported in literature (e.g.,~\cite{Heldring2021, tetzner_2024, Laviada2009, Heldring2014, Heldring2015}). In these cases, the ACAPP terminates prematurely when the error estimate falls below the implemented threshold, even though the actual error remains above the threshold. Hence, this phenomenon is also known as \emph{early convergence}.

We cannot resort to the tempting solution of just assembling the matrix block $A$ entirely. This would go against the core objective of matrix compression: reducing memory consumption. Therefore, most proposals to alleviate early-convergence issues implement more elaborate pivoting strategies and termination criteria. One of the main idea behind these extensions is to calculate more matrix elements than strictly necessary in ACAPP and use this information to improve error estimation and the pivot selection. This should be performed with care, because full pivoting yields quadratic algorithmic complexity~\cite {bebendorf2003adaptive}, whereas log-linear complexity can only be achieved with partial-pivoting strategies. Therefore, we propose only algorithmic extensions that do not degrade the computational complexity.

\subsection{Diagonal convergence criterion}

One of the core reasons for early convergence is that the termination criterion checks the compression error only on a part of the matrix block. As explained, calculating entire matrix blocks to verify convergence worsens the algorithm's complexity to quadratic. However, one could compute a small amount of additional matrix elements without increasing the overall log-linear complexity. Specifically, given a block of size $m \times n$, we can calculate $d$ additional matrix elements, where $d$ scales linearly with $m$ and $n$, and maintain the log-linear complexity of ACAPP. Those $d$ matrix elements are then used by the error estimator to improve the convergence criteria. Different approaches to choose which $d$ elements to calculate have been described in, e.g., \cite{Heldring2021} and \cite{tetzner_2024}.

Before proposing our approach to choosing $d$ additional matrix elements, let us explain how they are used in the error estimator. Let us consider a vector $\mathbf{e}_0=\{e_0^0,\dots,e_0^{d-1}\}$ with the extra entries $e_0^s = A[i_s,j_s]$ for $s\in\{0,\dots,d-1\}$. Here $(i_s,j_s)$ are the row index and column index of entry $e_0^s$, respectively. Then, we can update this vector in iteration $\ell$ of the algorithm with
\begin{equation}
    \mathbf{e}_{\ell}[s] := \mathbf{e}_{{\ell}-1}[s] - \mathbf{u}_{\ell}[i_s]\mathbf{v}_{\ell}[j_s],\quad \text{for }s\in\{0,\dots, d-1\}.
    \label{eq:update_error_extra_pivots}
\end{equation}
Using the notation $|\mathbf{e}_{\ell}^2| := \{ |(e_{\ell}^0)^2|,\dots,|(e_{\ell}^{d-1})^2| \}$, we define the additional convergence criterion as
\begin{equation}
    \sqrt{\operatorname{mean}\left(\left|\mathbf{e}_{\ell}^2\right|\right) m n} \leqslant \varepsilon\left\|R_{\ell}^{m \times n}\right\|_{\mathrm{F}},
    \label{eq:rscc}
\end{equation}
where the Frobenius norm of $R_{\ell}$ can be calculated with Equation~\eqref{eq:fronorm_approximation}. 

Selecting these additional $d$ matrix elements in each block is a design choice. For example, the random-sampling convergence criterion (RSCC) proposed in \cite{Heldring2021} uses a fixed value of $d$ and a random selection of matrix elements within each block. In~\cite{tetzner_2024}, they choose $d=m+n$ random elements to solve early-convergence problems for magnetic field integral equations while maintaining log-linear complexity.

\begin{figure}[htbp]
    \begin{subfigure}[b]{0.35\linewidth}
        \centering
        \includegraphics[height=4.1cm]{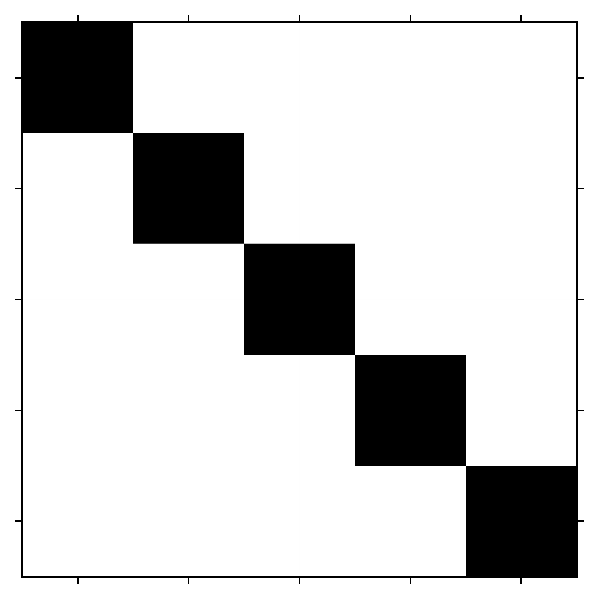}
        \subcaption[$5\times5$ matrix]{$5\times5$ matrix.}
    \end{subfigure}
    \hfill
    \begin{subfigure}[b]{0.28\linewidth}
        \centering
        \includegraphics[height=4.1cm]{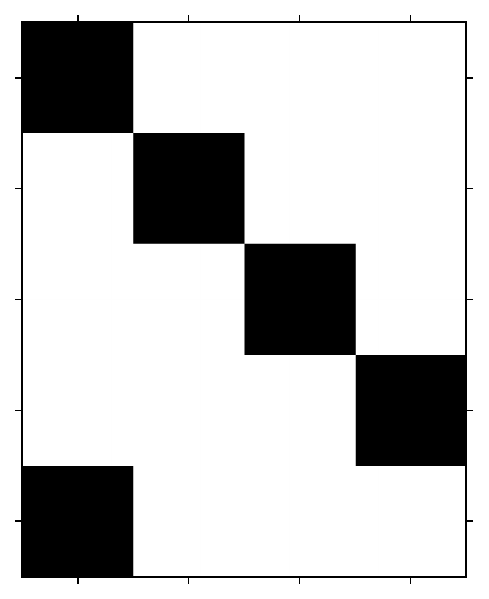}
        \subcaption[$5\times4$ matrix]{$5\times4$ matrix.}
    \end{subfigure}
    \hfill
    \begin{subfigure}[b]{0.23\linewidth}
        \centering
        \includegraphics[height=4.1cm]{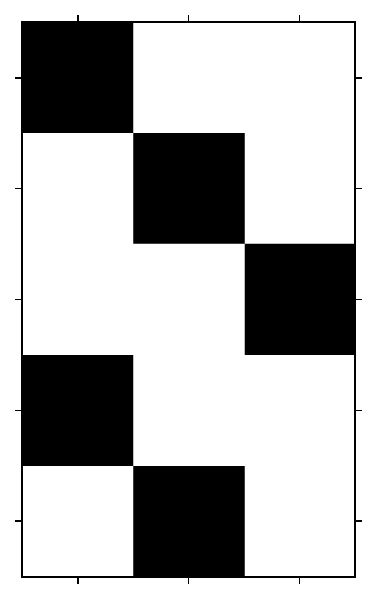}
        \subcaption[$5\times3$ matrix]{$5\times3$ matrix.}
    \end{subfigure}
    \begin{subfigure}[b]{1.0\linewidth}
        \centering
        \includegraphics[height=4.1cm]{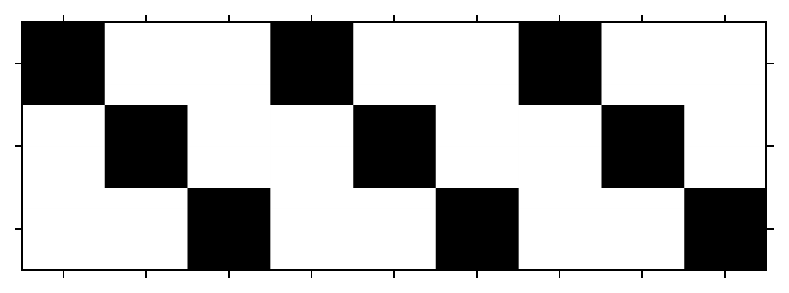}
        \subcaption[$3\times9$ matrix]{$3\times9$ matrix.}
    \end{subfigure}
    \caption{Examples of extended diagonals for different matrices. The coordinates of the extended diagonal are calculated using modulo operations with the number of rows and columns of the corresponding rectangular matrix.}
    \label{fig:extended_diagonals}
\end{figure}

In our implementation, we decided to use the extended main diagonal of the matrix as additional elements $\mathbf{e}_0$; see Figure~\ref{fig:extended_diagonals}. We call this the Diagonal Convergence Criterion (DCC). Notice that this involves calculating $d=\max(m,n)$ extra entries for each matrix block of size $m \times n$, thus avoiding an increase in storage complexity. Using a deterministic choice of elements (rather than random selections) improves robustness, consistency, and reproducibility of the compression algorithm. Furthermore, computational evidence indicates that early-convergence problems are associated with specific matrix structures. Examples of compression error patterns we have seen include periodic subblocks and groups of rows or columns. Moreover, the benefit of using the extended main diagonal is that information from all rows and columns are included, avoiding common error patterns.

We propose combining three convergence criteria: Equations~\eqref{eq:final_inequality} and~\eqref{eq:euclidean_relative_error} for two common error estimators, and Equation~\eqref{eq:rscc} for the proposed DCC. Our combined convergence criterion is thus defined as
\begin{equation}
    \max\left\{
    \frac{\|\mathbf{u}_{\ell}\|_2\|\mathbf{v}_{\ell}\|_2}{\|\mathbf{u}_1\|_2\|\mathbf{v}_1\|_2},
    \frac{\|\mathbf{u}_{\ell}\|_2\|\mathbf{v}_{\ell}\|_2}{\|R_{\ell}\|_{\mathrm{F}}}, 
    \frac{\sqrt{\operatorname{mean}\left(\left|\mathbf{e}_{\ell}^2\right|\right) m n}}{\|R_{\ell}\|_{\mathrm{F}}}
    \right\}
    \leqslant \varepsilon.
    \label{eq:combined_convergence_criterion}
\end{equation}
Algorithms that include this combined convergence criterion will have the suffix DCC, for example, ACAPP~DCC and ACAPP+~DCC.

\FloatBarrier

\subsection{Diagonal pivots} 

The primary reason for the $d$ additional matrix elements is to improve the error estimator. However, as they are already calculated, they can also be used in the pivoting strategy without incurring memory overhead. This idea was applied to the RSCC in \cite{tetzner_2024}, and we apply it to the DCC. In more detail, we detect the case in which the standard convergence criteria from Section~\ref{sec:error-estimation} are satisfied, but the DCC in Equation~\eqref{eq:rscc} is not yet reached in iteration $\ell$. This situation suggests that the default pivoting strategy is deficient and significant parts of the matrix are left untreated. A fundamentally different pivot selection with a broader perspective must be adopted to address unexplored regions of the matrix block that affect compression accuracy. Here the $d$ matrix elements of the DCC are used as candidates for the next iteration pivot ($i_\star$ and $j_\text{ref}$ for the ACAPP and ACAPP+, respectively), selecting the entry in $\mathbf{e}_{\ell}$ with the largest magnitude as the next pivot. Note that changing the value of $j_\text{ref}$ involves updating $\mathbf{u}^\text{ref}$, $i_\text{ref}$ and $\mathbf{v}^\text{ref}$ before the next iteration. Additionally, these pivot candidates can also be used in Line~\ref{alg_line:aca+_random_reference} of Algorithm~\ref{alg:aca_plus} instead of selecting a random reference column, removing the random component of the algorithm.

Importantly, as the DCC calculates the extended diagonal, the $d$ additional matrix elements cover all rows and columns. This guarantees that there is always a pivot in $\mathbf{e}_{\ell}$ that has not yet been identified by the algorithm and therefore has potentially valuable information for the compression. In fact, in the extreme case where all DCC pivots are used, the entire matrix is covered, and the decomposition is exact. Providing the diagonal elements as optional pivots thus improves the reliability and robustness of the compression algorithm.

\subsection{Extended admissibility condition}

\begin{figure}[htbp]
    \centering
    \includegraphics[width=.9\linewidth]{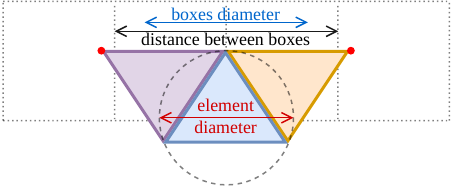}
    \caption{Visualization of a scenario where the additional condition~\eqref{eq:admissibility2} is necessary to avoid admissible boxes with near-singular interactions. Specifically, the two red vertices belong to adjacent triangular mesh elements, for which special singularity-aware quadrature is employed. At the same time, they belong to two octree boxes that satisfy the standard admissibility condition~\eqref{eq:admissibility}.
    }
    \label{fig:adjacent_elements}
\end{figure}

The integral operators depend on Green's function~\eqref{eq:green_function_helmholtz}, which has a singularity at zero distances. Therefore, integrands are singular when evaluating self-interactions and nearly-singular for interactions between neighboring mesh elements. Importantly, even though the integrand may be singular, the integrals are well defined. Still, special singularity-aware quadrature schemes must be used to achieve accurate numerical evaluations for self and near-interactions~\cite{sauter2010boundary}. In practice, our P1-BEM implementation selects a singular or regular quadrature scheme based on the mesh topology, with all triangle pairs that share at least one vertex being handled by singularity-aware quadrature. This means that, in rare cases, near interactions may occur within admissible blocks. See Figure~\ref{fig:adjacent_elements} for an example of P1 elements on a triangular surface mesh. Calculating nearly-singular integrals in the ACAPP is undesirable for two reasons. First, checking for near singularity during the ACAPP increases computation time and prevents vectorization. Second, including all near interactions in the compressed matrix is essential to achieving overall accuracy, but they may not be selected by the ACAPP pivots when in admissible blocks.

To achieve robust matrix compression, we decided that all singularity-aware quadrature, also for near interactions, must be handled by dense matrix blocks, and never by the ACAPP. We achieved this robustness requirement by adjusting the admissibility condition. Specifically, the admissibility condition in Equation~\eqref{eq:admissibility} still allows for nearly-singular integrals between local test and basis functions that share a common edge or vertex in the mesh to end up in admissible blocks. To prevent this behavior, we also require
\begin{equation}
    \operatorname{dist}(Box_i, Box_j) > 2\ \text{max\_element\_diameter}
    \label{eq:admissibility2}
\end{equation}
for boxes to be admissible. We calculate the element diameter as the diameter of the circumscribed circle or sphere around a surface or volume grid element, respectively. The maximum element diameter is the largest among all elements in the mesh. Notice that the multiplicative factor of 2 in condition~\eqref{eq:admissibility2} is necessary because two points of the mesh can be in two different boxes with a distance larger than the maximum element diameter, yet still belong to adjacent elements. An example of this situation is depicted in Figure~\ref{fig:adjacent_elements}.

\subsection{Sustained convergence check}

Another effect of the approximate error estimator is that the ACAPP's convergence is not always monotone. That is, while adding new vectors to the low-rank approximation always improves the exact compression error, it may actually worsen for the error estimator. This also means that the ACAPP may satisfy the tolerance at iteration $\ell$ but not at $\ell+1$ or subsequent iterations. Again, in these rare cases, the ACAPP converges too early.

We propose a sustained convergence check. Under this scheme, the convergence criterion must hold for $\alpha$ consecutive iterations before we terminate the ACAPP. This stronger termination criterion avoids cases in which a single iteration falls below the convergence criterion, but subsequent ones do not. In fact, it avoids even rarer cases where $\alpha-1$ consecutive iterations have converged while the next one has not. Hence, the ACAPP terminates when all iterations from $\ell$ to $\ell+\alpha-1$ satisfy the convergence criterion. Even though a low-rank approximation of rank $\ell+\alpha-1$ is available, we choose to eliminate the last vectors and store the low-rank approximation of the first $\ell$ vectors. This avoids memory overhead.

Similar ideas can be found in the literature. For example, multiple restarts are used in~\cite{Laviada2009}. In that case, additional iterations are also performed when the algorithm reaches the desired error. However, they use random initial pivots for each restart. In contrast, our proposal is deterministic and does not require restarts.

Our computational experiments indicate that choices as small as $\alpha = 2$ or $\alpha = 3$ are sufficient to achieve robust compression. We emphasize that although the calculation time is slightly longer, this does not incur memory overhead in the final compressed matrix, and the algorithm is more robust.

\section{Results}\label{sec:results}

This section presents numerical evidence that our proposed algorithmic extensions to the ACAPP matrix compression technique improve robustness while maintaining efficiency. We showcase several benchmarks that exhibit early-convergence problems and demonstrate that the ACAPP+~DCC variant remains robust across all experiments. A parametric study confirms the log-linear complexity of our implementation. Finally, we demonstrate the practical efficacy of matrix compression on a CT-derived human skull model.

\subsection{Computational settings}

The hierarchical matrix compression methodology is implemented in Python (v.~3.12.2), using the Numpy library (v.~1.26.4) \cite{numpy} for linear algebra, and the Numba library (v.~0.59.1) \cite{numba} for shared-memory parallel computing via the multithreading paradigm. The BIOs with P1 elements are calculated using the Bempp-cl library (v.~0.4.2) \cite{bempp-cl}. The VIOs and BIOs with P0 elements are programmed in Python, and have been verified against analytical solutions and finite element methods in~\cite{aballay2026nested}. All numerical experiments were performed with double-precision floating point numbers on a compute node with 32~cores shared across two Intel(R) Xeon(R) Silver 4216 2.10~GHz processors and a total of 2048~GB of RAM.

\subsection{Compressed matrix structure}
\label{sec:compressed_structure}

The hierarchical matrix compression technique approximates integral operators corresponding to admissible boxes in the octree with low-rank decompositions. Let us visualize the resulting compression structure to assess the effectiveness of our implementation. As a representative example, we consider the DL-BIO with P1 and P0 elements (denoted by P1-DL-BIO and P0-DL-BIO, respectively) on the surface of a unit sphere. The P1-DL-BIO's triangular surface mesh has 24,845 vertices, and the P0-DL-BIO's voxel mesh has 31,014 rectangular faces on the surface; see Figure~\ref{fig:sphere_meshes}. The nondimensional material parameters are provided in Table~\ref{tab:sphere_parameters}. The results are compared for two $\varepsilon$ values, $10^{-1}$ and $10^{-10}$. A maximum depth of 5 and minimum box sizes of 15 and 20 were taken for the P1-DL-BIO and P0-DL-BIO, respectively. For the sustained convergence check, an $\alpha$ value of 2 was used for both experiments.

\begin{table}[htbp]
\caption{Nondimensional parameters used in the benchmarks on a sphere with unit radius.}
\label{tab:sphere_parameters}
\centering
\begin{tabular}{lr}
\toprule
parameter & value \\
\midrule
density (exterior)        & 1     \\ 
speed of sound (exterior) & 1     \\ 
density (interior)        & 2     \\ 
speed of sound (interior) & 3     \\ 
frequency                 & 6.7   \\ 
points per wavelength     & 6     \\
\bottomrule
\end{tabular}
\end{table}

\begin{figure}[htbp]
    \begin{subfigure}[b]{.49\linewidth}
        \includegraphics[width=\linewidth]{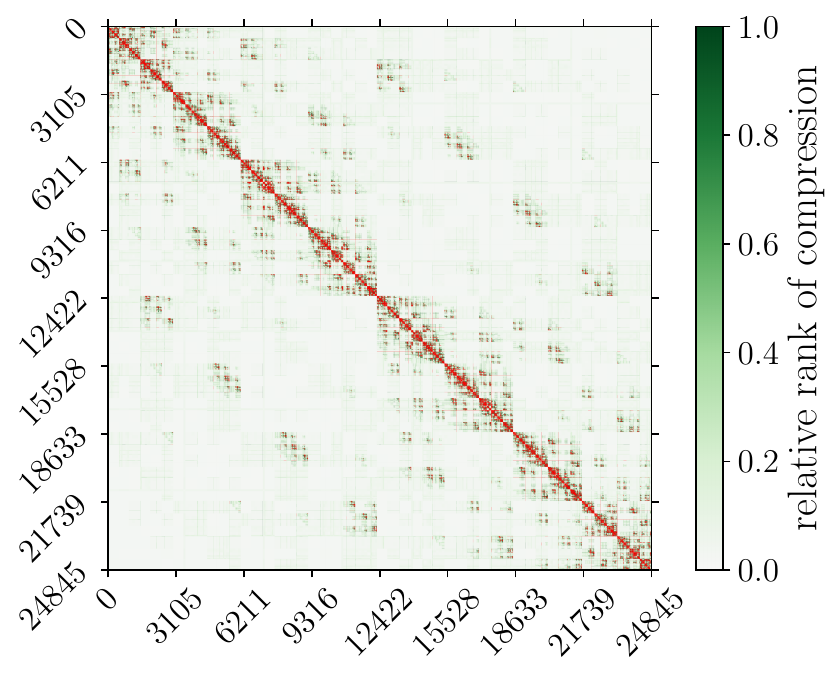}
        \subcaption{P1: $\varepsilon=10^{-1}$}
    \end{subfigure}
    \hfill
    \begin{subfigure}[b]{.49\linewidth}
        \includegraphics[width=\linewidth]{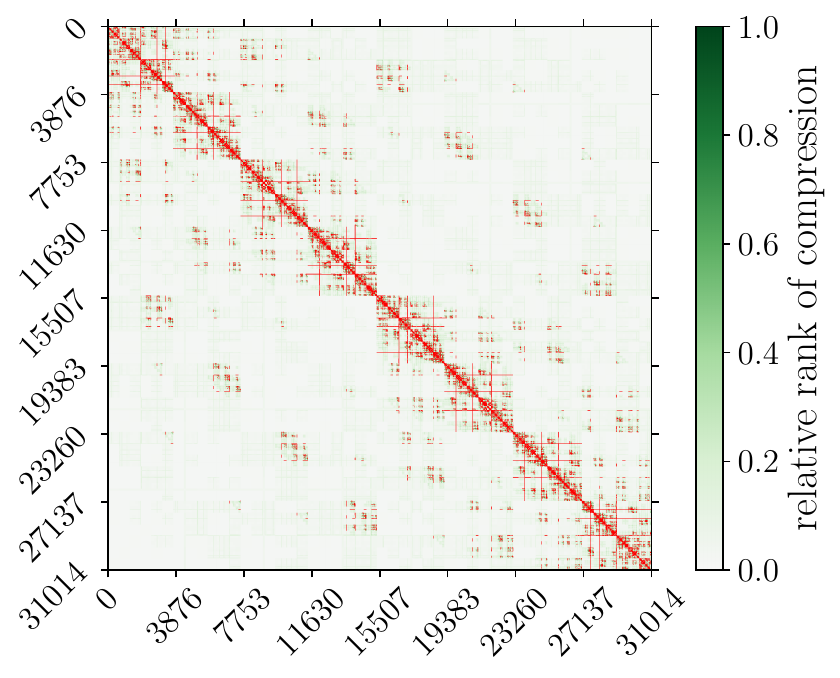}
        \subcaption{P0: $\varepsilon=10^{-1}$}
    \end{subfigure}

    \begin{subfigure}[b]{.49\linewidth}
        \includegraphics[width=\linewidth]{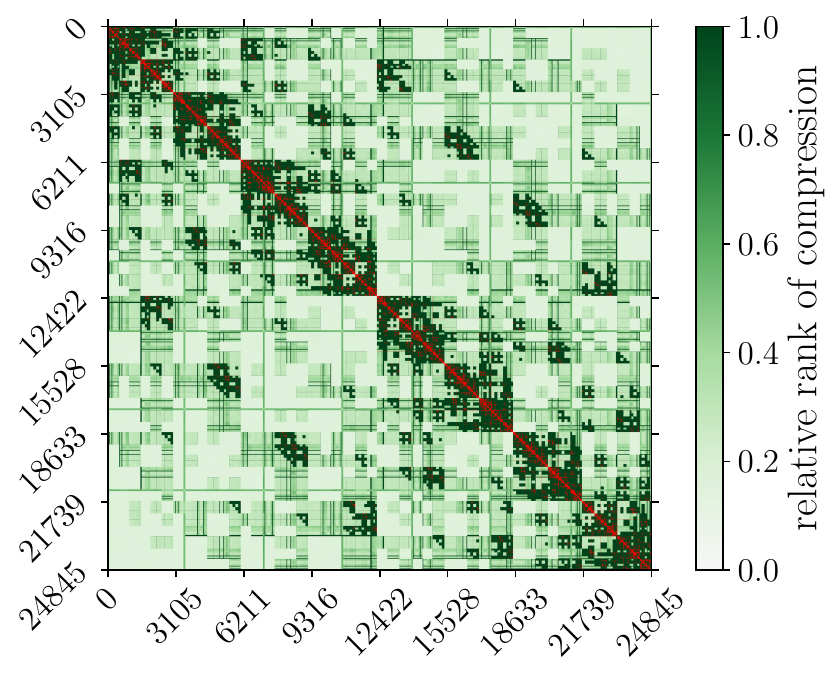}
        \subcaption{P1: $\varepsilon=10^{-10}$}
    \end{subfigure}
    \hfill
    \begin{subfigure}[b]{.49\linewidth}
        \includegraphics[width=\linewidth]{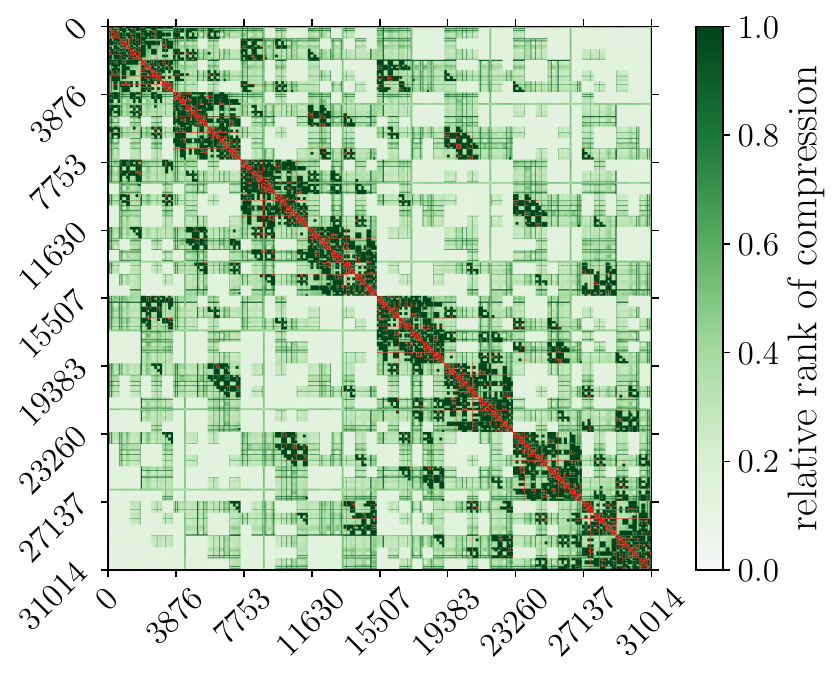}
        \subcaption{P0: $\varepsilon=10^{-10}$}
    \end{subfigure}
    
    \caption{Matrix compression patterns for the ACAPP+~DCC for the DL-BIOs on the surface of a unit sphere and $\varepsilon\in\{10^{-1},10^{-10}\}$. Panels (a) and (c) correspond to the P1 results, and panels (b) and (d) to the P0 results. The rows and columns of the matrix were sorted according to a depth-first search of the leaves of the corresponding octree. The red pixels indicate non-admissible interactions.}
    \label{fig:sphere_imshows}
\end{figure}

Figure~\ref{fig:sphere_imshows} displays the hierarchical matrix structure resulting from the compression. Red pixels denote non-admissible blocks that are stored in their full dense form. Green areas depict the compressed, admissible blocks. The shade of the green color scale represents the achieved relative rank of compression, which is calculated by dividing the compression rank by the full dense rank of the corresponding matrix block. In other words, the darker shades indicate a smaller reduction in rank (lower compression effectiveness), whereas lighter shades correspond to a larger rank reduction (higher compression effectiveness).

These visualizations confirm the hierarchical block structure generated by the octree subdivision. As observed, admissible interactions constitute the vast majority of the matrix elements, whereas inadmissible blocks are primarily restricted to locations close to the main diagonal and small clusters scattered throughout the rest of the matrix. Regarding the error thresholds, there is a clear connection between the prescribed tolerance and the resulting block ranks: as $\varepsilon$ decreases, the ranks of the compressed blocks increase. Consequently, stricter accuracy requirements inherently lead to more memory consumption and lower compression efficiency.

Additionally, the two discretization strategies of triangular P1 versus rectangular P0 elements yield consistent compression results, with only slight deviations at certain matrix locations. The observed similarities are expected as both formulations employ the same integral operator. Furthermore, the hierarchical octrees used for compression are largely identical at their upper levels, with the main differences appearing only at the leaf boxes that represent the smallest mesh partitions.

\subsection{Compression error}

The objective of matrix compression is to reduce memory consumption while controlling the error incurred in a matrix-vector multiplication. However, the ACAPP for low-rank approximations does not guarantee error bounds as it relies on error estimators based on partial matrix information. Hence, let us numerically examine the compression error for different values of the compression tolerance, i.e., $\varepsilon \in \{10^{-1}, 10^{-2}, \dots, 10^{-10}\}$. We measure the relative error as
\begin{equation}
    E^\varepsilon_\text{rel}(\mathbf{v})=\frac{\|A\mathbf{v}-A_\varepsilon \mathbf{v}\|_2}{\|A\mathbf{v}\|_2},
    \label{eq:error_matvec}
\end{equation}
where $A$ is the original dense matrix of the numerical system, $A_\varepsilon$ is the compressed matrix for an error threshold of $\varepsilon$, and $\mathbf{v}$ is a random vector. In our experience, the error measure was not noticeably affected by using a different seed for the random number generator or by averaging over multiple random vectors. Hence, we present the results for a single random vector~$\mathbf{v}$.

\begin{figure}[htbp]
    \begin{subfigure}[b]{0.497\linewidth}
        \centering
        \includegraphics[width=\linewidth]{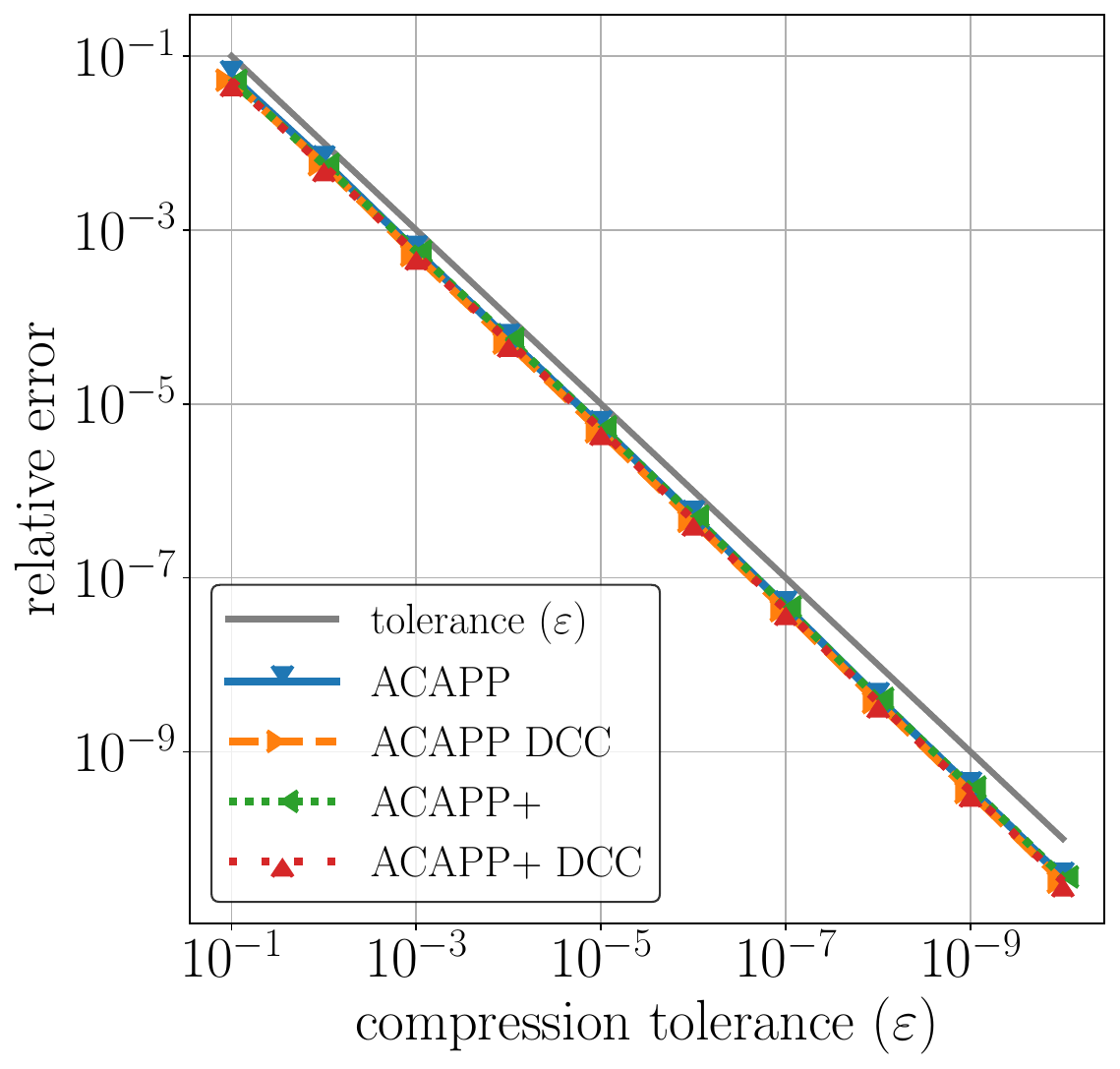}
        \subcaption{P1 discretization.}
    \end{subfigure}
    \hfill
    \begin{subfigure}[b]{0.497\linewidth}
        \centering
        \includegraphics[width=\linewidth]{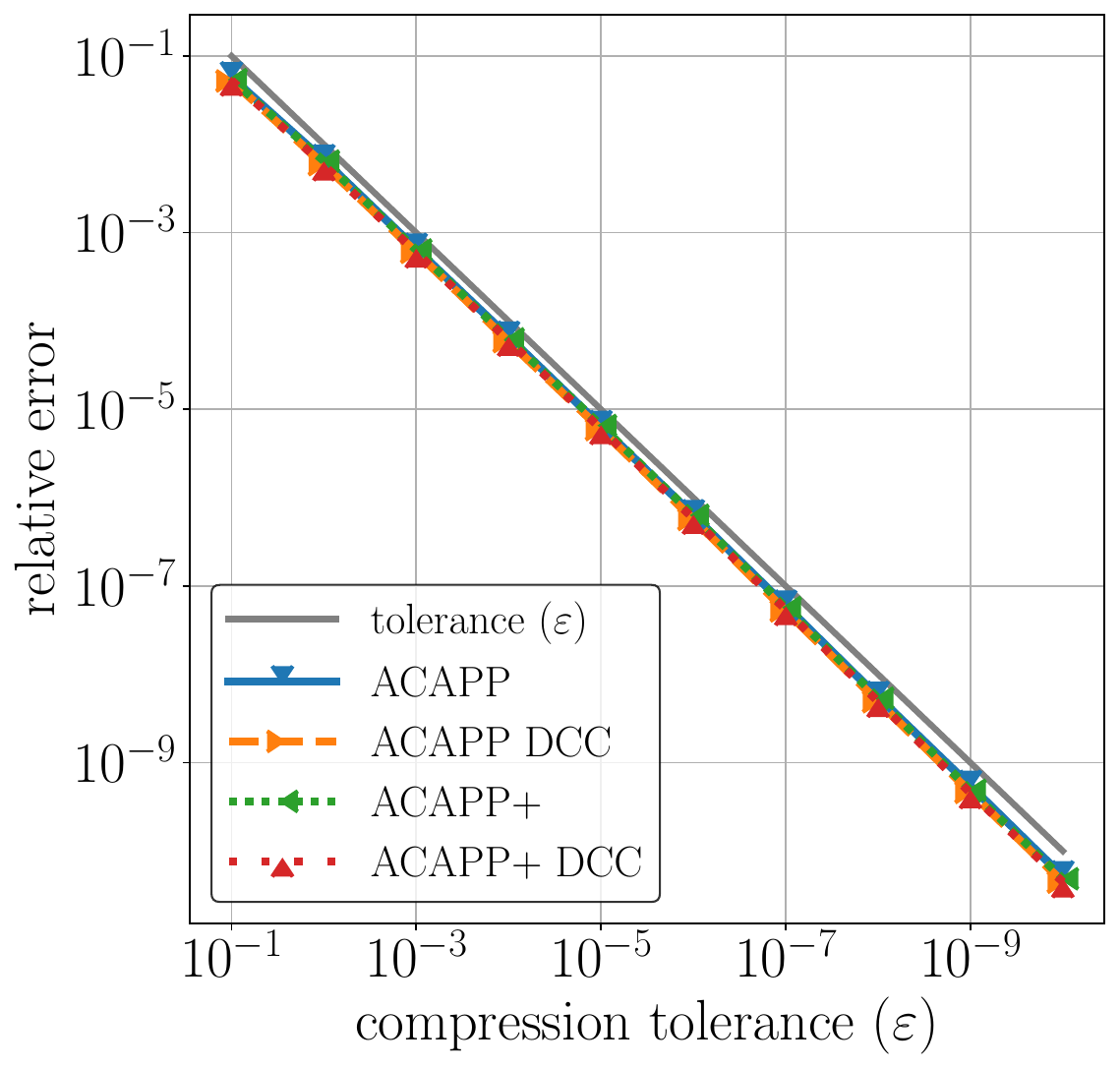}
        \subcaption{P0 discretization.}
    \end{subfigure}
    \caption{Relative errors~\eqref{eq:error_matvec} obtained for the compression of the DL-BIO on the surface of a unit sphere, with material parameters as in Table~\ref{tab:sphere_parameters}, and compression tolerances $\varepsilon\in\{10^{-1},10^{-2},\dots,10^{-10}\}$. The gray line represents the specified error tolerance for each compression.}
    \label{fig:sphere_relative_errors}
\end{figure}

Figure~\ref{fig:sphere_relative_errors} presents the relative errors for the DL-BIOs on the surface of the sphere, computed according to Equation~\eqref{eq:error_matvec}. In both cases, the observed errors of all four ACAPP versions remain consistently below the prescribed tolerances. These results verify the efficacy of the algorithm's tolerance-bound compression error. In our experience, the compression error almost always remains below the prescribed tolerance. However, there are specific cases in which some ACAPP versions fail to compress matrix blocks sufficiently.

\subsection{Early convergence}
\label{sec:early-convergence-results}

Despite the previous results showing that all four versions of the compression algorithm perform well on the P1-DL-BIO and P0-DL-BIO sphere benchmarks, there are rare cases in which standard pivoting strategies fail due to early convergence. Let us consider the following two scenarios: 1) a P0-DL-BIO on the surface of a unit cube, and 2) a P0-SL-VIO on the interior of an ellipsoid with semi-axes of length 1, 1, and 5. In both cases, we take the material parameters from Table~\ref{tab:sphere_parameters}, but the frequencies are 20 and 3, respectively. A maximum tree depth of 5 was maintained for both cases, with a minimum box size of 20 for the P0-DL-BIO on the cube and 30 for the P0-SL-VIO on the ellipsoid. We used $\alpha = 3$ for the sustained convergence check.

\begin{figure}[htbp]
    \begin{subfigure}[b]{0.497\linewidth}
        \centering
        \includegraphics[width=\linewidth]{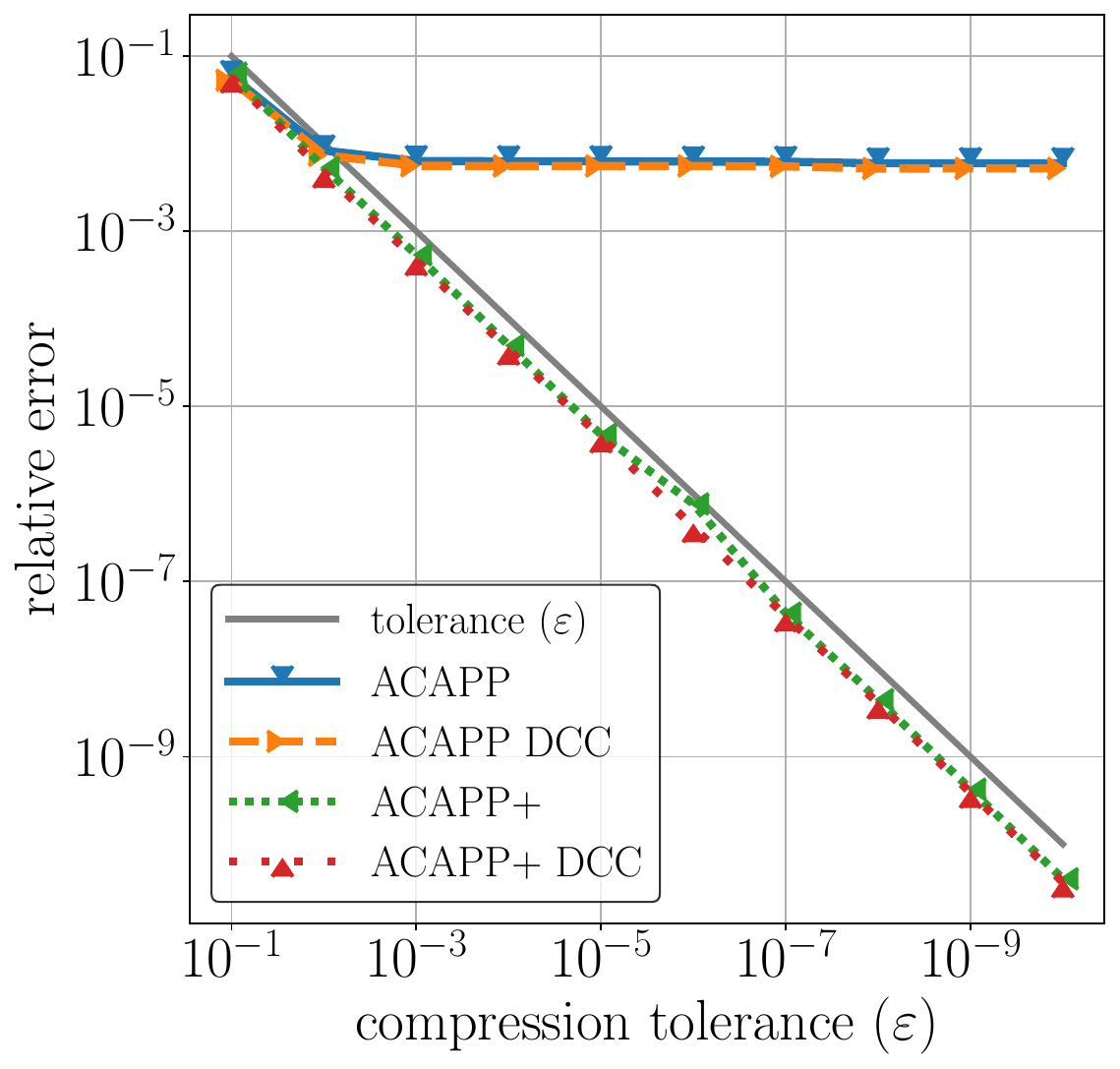}
        \subcaption{P0-DL-BIO on a cube surface.}
        \label{fig:early_convergence_relative_errors_cube}
    \end{subfigure}
    \hfill
    \begin{subfigure}[b]{0.497\linewidth}
        \centering
        \includegraphics[width=\linewidth]{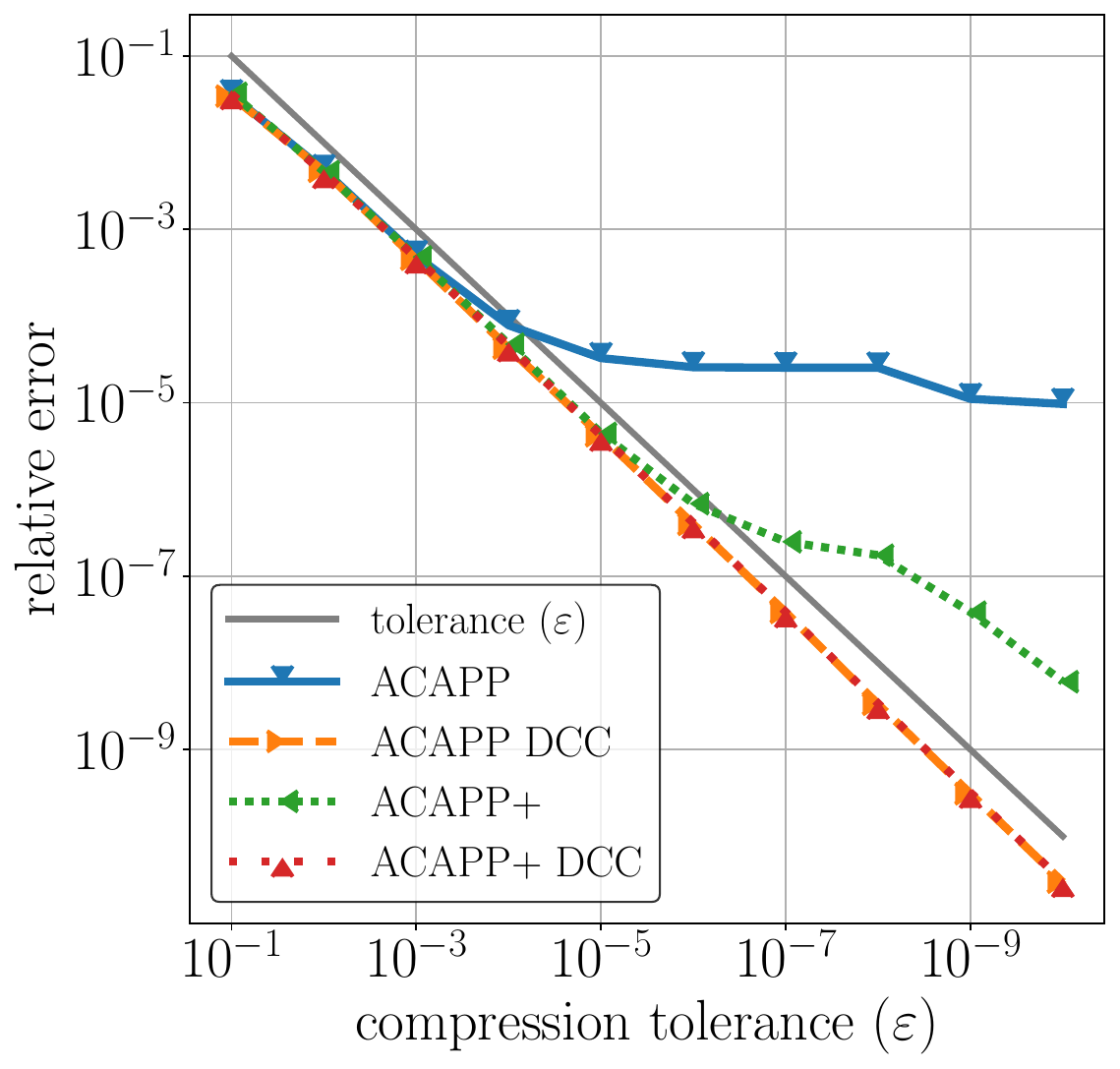}
        \subcaption{P0-SL-VIO on an ellipsoid.}
        \label{fig:early_convergence_relative_errors_ellipsoid}
    \end{subfigure}
    \caption{Relative errors obtained for the compression of the P0-DL-BIO on the surface of a unit cube and the P0-SL-VIO on the interior of an ellipsoid with semi-axes of lengths (1,1,5) for $\varepsilon\in\{10^{-1},10^{-2},\dots,10^{-10}\}$ in a matrix-vector operation with a random vector. The gray line represents the specified error tolerance for each compression.}
    \label{fig:early_convergence_relative_errors}
\end{figure}

Figure~\ref{fig:early_convergence_relative_errors} clearly illustrates the impact of early convergence that was discussed in Section~\ref{sec:error-estimation}. Essentially, the ACAPP algorithm significantly underestimates the compression error of some matrix blocks. This causes the ACAPP algorithm to prematurely terminate its iterative creation of low-rank approximations, well before reaching the target tolerance. Importantly, while this phenomenon does not occur in all blocks, the poor approximation of even a single block may cause significant errors for the global matrix-vector product. Strikingly, the ACAPP+~DCC variant does satisfy the compression tolerance in both benchmarks.

\begin{figure}[htbp]
    \begin{subfigure}[b]{0.495\linewidth}
        \centering
        \includegraphics[width=\linewidth]{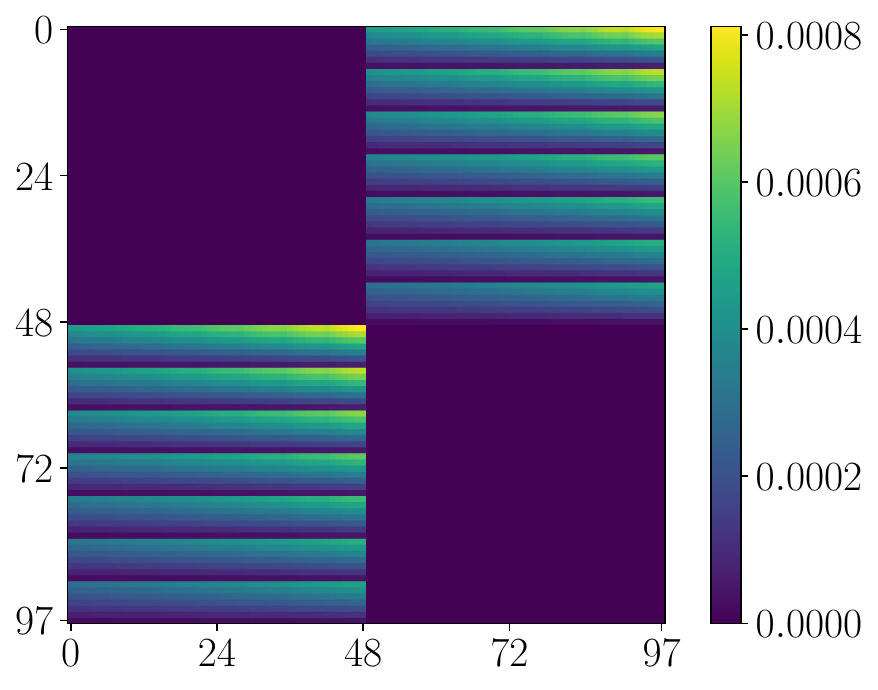}
        \subcaption{Original matrix block.}
    \end{subfigure}
    \hfill
    \begin{subfigure}[b]{0.495\linewidth}
        \centering
        \includegraphics[width=\linewidth]{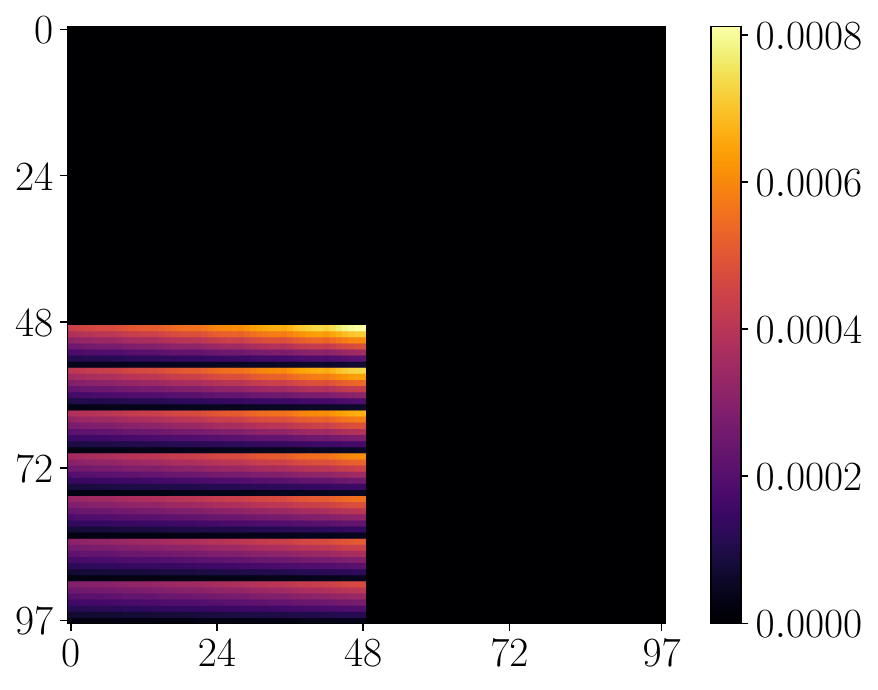}
        \subcaption{Elementwise compression error.}
    \end{subfigure}
    
    \caption{An example of a matrix block that presents early convergence in the standard ACAPP variant, for a P0-DL-BIO on the surface of a unit cube and a tolerance of $\varepsilon=10^{-4}$. The color of the matrix elements represents the absolute value of the original block (left) and the elementwise absolute difference between the original block and the compression (right).}
    \label{fig:bad_compression_examples_cube}
\end{figure}

One of the most prominent instances of early convergence we encountered is a structure specific to the DL-BIO on a cube, which is consistent with the literature on double-layer operators in domains with edges~\cite{grasedyck2005}. As is clear from Figure~\ref{fig:early_convergence_relative_errors}, at tolerance values below $\varepsilon = 10^{-3}$, compression stalls for the ACAPP and ACAPP~DCC algorithms. To understand the reason behind this behavior, let us investigate the structure of one of such matrix blocks that causes the early-convergence problem; see Figure~\ref{fig:bad_compression_examples_cube}. The matrix blocks are characterized by two non-zero sub-blocks separated by regions of zero elements. These zero sub-blocks are common in double-layer operators (see Equation~\eqref{eq:dl-bio}), where dot products between perpendicular normal vectors from adjacent faces of the cube yield zero-valued matrix elements. Consequently, the standard pivoting strategy in ACAPP and ACAPP~DCC extracts information from only one nonzero sub-block. The ACAPP's error estimation will therefore consider only that nonzero sub-block. Hence, convergence will be achieved without considering the other nonzero sub-block at any stage of the algorithm. 
Indeed, the elementwise error displayed in Figure~\ref{fig:bad_compression_examples_cube} confirms that one of the nonzero sub-blocks was approximated within the specified tolerance. 
At the same time, the other nonzero sub-block has a relative elementwise error of 100\% since it was left untreated in the low-rank decomposition. This problem cannot be solved by lowering the compression tolerance. In contrast, the two ACA variants employing the improved pivoting strategy (ACAPP+ and ACAPP+~DCC) overcome this problem by simultaneously evaluating information from both sub-blocks during the pivoting strategy in each iteration of the low-rank approximation. They achieve a global representation of the matrix block in the compression error estimate and compress within the specified tolerance.

\begin{figure}[htbp]
    \begin{subfigure}[b]{0.52\linewidth}
        \centering
        \includegraphics[width=\linewidth]{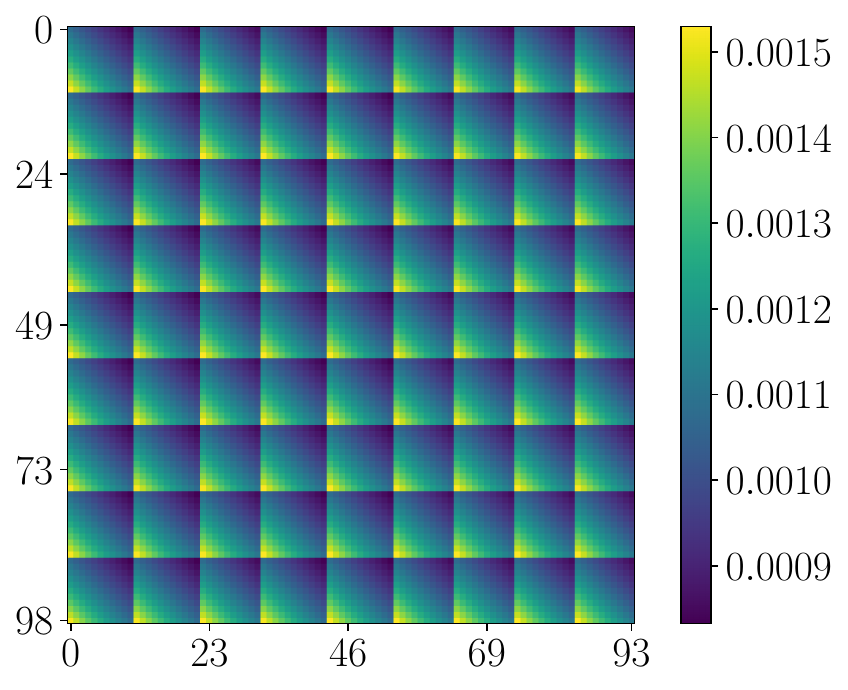}
        \subcaption{Original matrix block.}
    \end{subfigure}
    \hfill
    \begin{subfigure}[b]{0.47\linewidth}
        \centering
        \includegraphics[width=\linewidth]{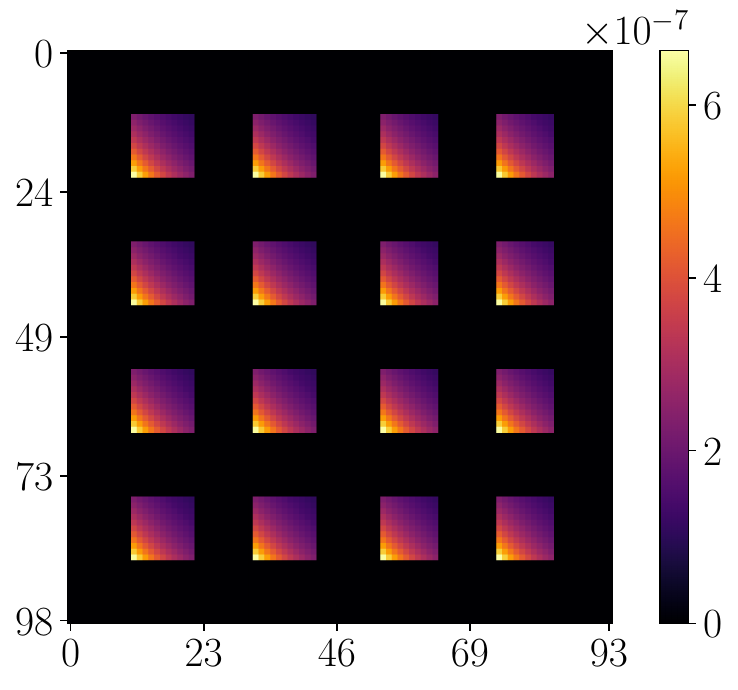}
        \subcaption{Elementwise compression error.}
    \end{subfigure}
    
    \caption{An example of a matrix block that presents early convergence in the standard ACAPP variant, for a P0-SL-VIO on the interior of an ellipsoid with semi-axes of lengths (1,1,5) and a tolerance of $\varepsilon=10^{-6}$. The color of the matrix elements represents the absolute value of the original block (left) and the elementwise absolute difference between the original block and the compression (right).}
    \label{fig:bad_compression_examples_ellipsoid}
\end{figure}

Another instance of early convergence can be observed in the compression of the SL-VIO for an ellipsoid. In this scenario, early-convergence errors arise from fully dense blocks that contain repeating sub-structures throughout the matrix, like the one shown in Figure~\ref{fig:bad_compression_examples_ellipsoid}. This specific arrangement of matrix information can lead compression algorithms to overlook some of these sub-structures during their execution. Consequently, integrating the DCC technique with the extended diagonal enables the algorithms to overcome this issue by incorporating specifically placed pivots designed to identify such information losses. Meanwhile, algorithms without DCC continue to exhibit early-convergence problems, as illustrated in Figure~\ref{fig:early_convergence_relative_errors_ellipsoid}.

\subsection{Compression rates}

Thus far, we have observed that our low-rank approximations effectively approximate the dense matrices of the BIO and VIO systems. The ACAPP+~DCC algorithm demonstrates particularly high fidelity, even in scenarios where other variants fail due to early convergence. Now, let us evaluate the compression rates achieved as $\varepsilon$ decreases, quantifying the memory savings gained by hierarchical matrix compression.

\begin{figure}[htbp]
    \centering
    \includegraphics[width=\linewidth]{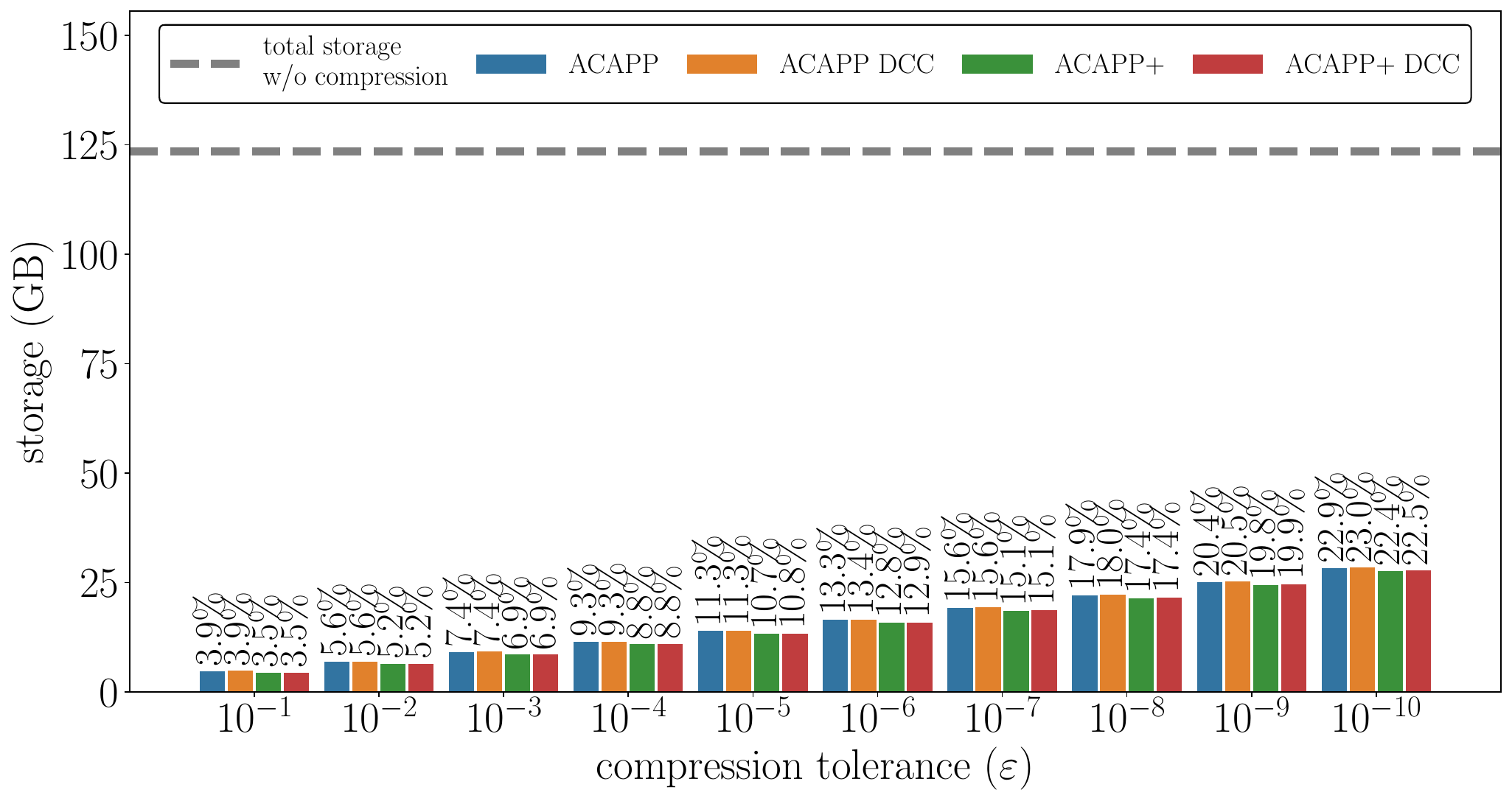}
    \caption{Memory storage and compression rates of the ACA algorithms, for a P0-DL-BIO on the surface of a unit cube. The gray line depicts the storage of the original dense matrix.}
    \label{fig:cube_final_storages}
\end{figure}

Figure~\ref{fig:cube_final_storages} displays the memory consumption and compression rates for the four ACA variants and the same double-layer operator of the unit cube presented in Section~\ref{sec:early-convergence-results}. Here, the compression rate is the memory consumption of the compressed matrix divided by the memory of the dense matrix. Notably, the compression rates remain similar across all four ACA variants. This indicates that the enhanced pivoting and convergence strategies used in the ACAPP+~DCC variant alleviate early-convergence problems without increasing memory consumption. In other words, the convergence problem was successfully addressed without compromising the compression rates. Furthermore, as expected, the required storage increases as the compression tolerance decreases.

For this scenario, the original dense matrix uses around 123~GB. At an $\varepsilon$ of $10^{-10}$, the compressed version requires only 28~GB, achieving a storage saving of 77.5\% even with such a strict tolerance.

\subsection{Storage complexity of matrix compression}

The discretized VIO and BIO are dense matrices that require $\mathcal{O}(n^2)$ memory, for $n$ the number of DOFs. Hierarchical matrix compression promises to reduce the storage complexity to $\mathcal{O}(n \log(n))$, which is essential to perform large-scale simulations. Let us study the storage complexity of the proposed algorithms for the VIOs. We use a cube with side lengths of 1~cm and an iterative mesh refinement of $2^j$~voxels per side, with $j\in\{2,\dots,8\}$. The materials are considered homogeneous with parameters resembling water and bone; see Table~\ref{tab:complexity_study_properties}.

\begin{table}[htbp]
    \centering
    \caption{Material parameters for the storage complexity study on a cube with 1~cm side lengths.}
    \label{tab:complexity_study_properties}
    \begin{tabular}{lrr}
    \toprule
    parameter & exterior & interior \\
    \midrule
    density (kg/m$^3$)             & 1000                         & 1787                         \\
    speed of sound (m/s)           & 1500                         & 2646                         \\
    \bottomrule
    \end{tabular}
\end{table}

We perform two experiments that handle the frequency and mesh resolution differently. The first experiment uses a constant frequency of 352.8~kHz, at which the wavelength is 7.5~mm in the interior. Hence, each mesh refinement increases the number of elements per wavelength. The second experiment follows the approach of fixing the number of voxels per wavelength while increasing the frequency. Precisely, we use 6 voxels per interior wavelength. The benchmark considers the compression of the matrix $B = T_{\text{P0}} - S_{\text{P0}}$, where $T_{\text{P0}}$ denotes the discrete P0-ADL-VIO and $S_{\text{P0}}$ the discrete P0-SL-VIO. This specific matrix is the largest block of a standard VSIE system~\cite{ultrasonics2026, aballay2026nested}.

We measure memory algebraically, not with software that tracks RAM usage. That is, given the number of DOFs $n$, storing the dense matrix requires $n^2$ complex-valued floating-point numbers of 128~bits. The storage of the compressed matrix is calculated from the ranks achieved in each matrix block. We perform the ACA algorithm on each block and store the achieved rank along with the block's dimension. This provides the number of matrix elements for the low-rank approximation of each block. These are then aggregated over the entire octree. This allows us to flush the memory after calculating a hierarchical block, and we do not need to store the entire compressed matrix. The benchmarks are performed for tolerance levels of $\varepsilon = 10^{-1}$ and $10^{-4}$, using the ACAPP+~DCC algorithm. In these experiments, we did not specify a maximum tree depth to fully exploit the compression potential. However, a minimum box size of 10 was maintained to prevent the proliferation of tiny blocks, which would otherwise increase the computational overhead of the tree without providing significant storage gains. No additional iterations were used for the sustained convergence check ($\alpha = 1$).

\begin{figure}[htbp]
    \centering
    \includegraphics[width=\linewidth]{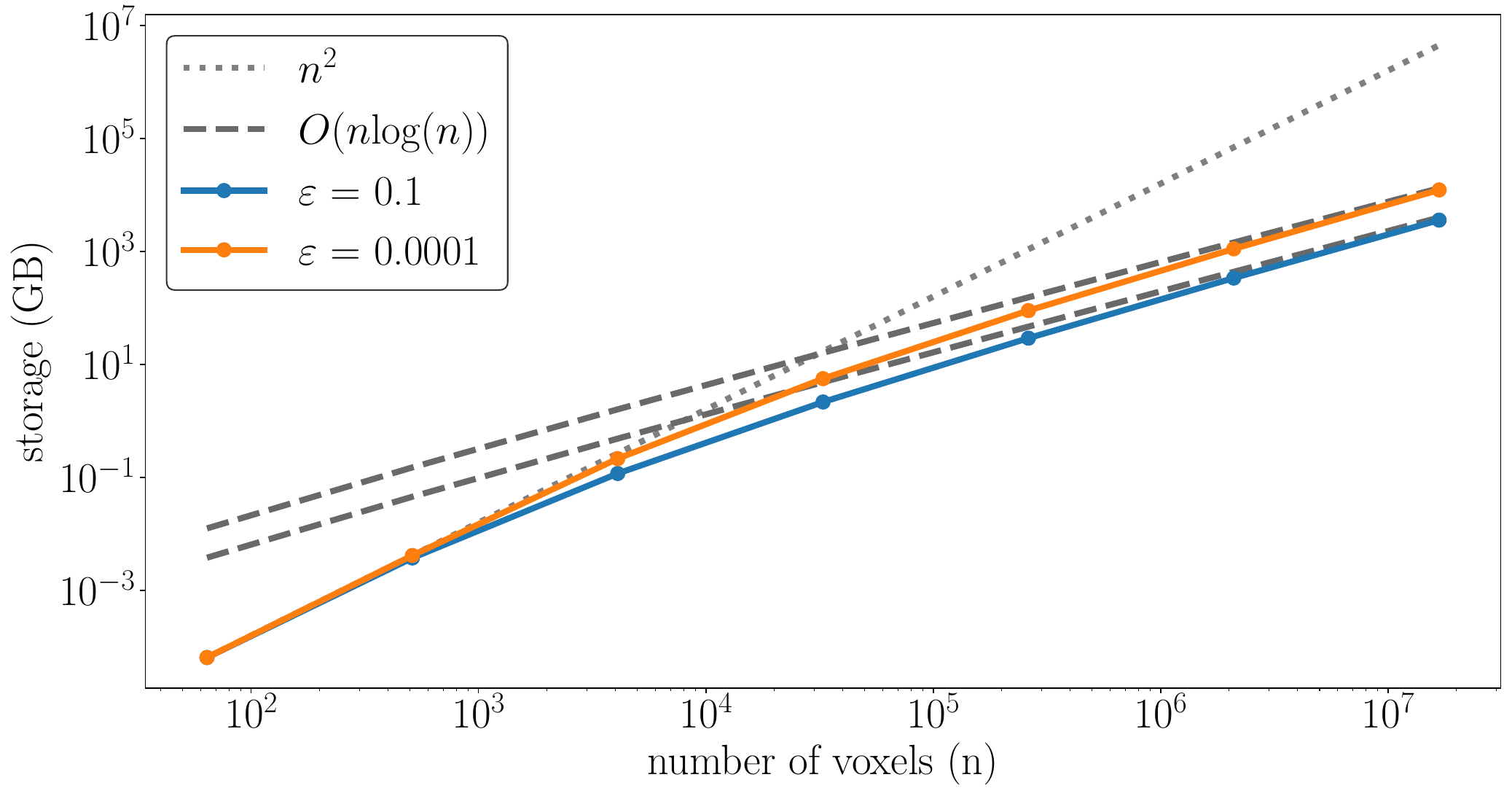}
    \caption{The memory requirements of the compression of matrix $B$ with tolerance $\varepsilon$ while using a fixed frequency at multiple grids. The dashed grey lines depict an $\mathcal{O}(n\operatorname{log}(n))$ complexity with different constants.}
    \label{fig:complexity_Fconst}
\end{figure}

\begin{table}[htbp]
\centering
\caption{The memory requirements of the compression of matrix $B$ with tolerance $\varepsilon$ while using a fixed frequency at multiple grids. The cubic domain has sides measuring 1~cm. Here, `vpw' is voxels per wavelength.}
\label{tab:complexity_Fconst}
\begin{tabular}{rrr|rr|rr}
\multicolumn{1}{l}{\multirow{2}{*}{\begin{tabular}[c]{@{}l@{}}number of\\ voxels ($n$)\end{tabular}}} &
  \multicolumn{1}{l}{\multirow{2}{*}{\begin{tabular}[c]{@{}l@{}}dense matrix\\storage (GB)\end{tabular}}} &
  \multicolumn{1}{l|}{\multirow{2}{*}{\begin{tabular}[c]{@{}l@{}}vpw\end{tabular}}} &
  \multicolumn{2}{c|}{$\varepsilon = 0.1$} &
  \multicolumn{2}{c}{$\varepsilon = 0.0001$} \\
\multicolumn{1}{l}{} &
  \multicolumn{1}{l}{} &
  \multicolumn{1}{l|}{} &
  \multicolumn{1}{l}{storage (GB)} &
  \multicolumn{1}{l|}{rate (\%)} &
  \multicolumn{1}{l}{storage (GB)} &
  \multicolumn{1}{l}{rate (\%)} \\ \hline
64         & 0.0000655     & 3   & 0.0000655 & 100.0 & 0.0000655  & 100.0 \\
512        & 0.00419     & 6   & 0.00379 & 90.42 & 0.00419  & 100.0 \\
4,096      & 0.27         & 12  & 0.12     & 44.11 & 0.22      & 80.88 \\
32,768     & 17.18        & 24  & 2.18     & 12.7  & 5.66      & 32.92 \\
262,144    & 1,099.51     & 48  & 29.44    & 2.68  & 90.65     & 8.24  \\
2,097,152  & 70,368.74    & 96  & 341.06   & 0.48  & 1,120.25  & 1.59  \\
16,777,216 & 4,503,599.63 & 192 & 3,617.98 & 0.08  & 12,359.05 & 0.27 
\end{tabular}
\end{table}

\FloatBarrier

As illustrated in Figure~\ref{fig:complexity_Fconst} and Table~\ref{tab:complexity_Fconst}, compression rates improve as the problem size increases. Specifically, while the number of DOFs increases cubically, the octree depth grows only logarithmically, for example, reaching a depth of 8 in the largest $2^8$ voxels-per-side case. Hence, the number of admissible interactions per level is bounded by the octree structure and the admissibility condition. Similarly, the compression ranks of the admissible blocks also remain low in the case of constant frequency and fixed $\varepsilon$. These factors collectively confirm that the memory storage of the algorithm approaches the log-linear complexity of $\mathcal{O}(n\log(n))$ as the mesh is refined for this experiment.

\begin{figure}[htbp]
    \centering
    \includegraphics[width=\linewidth]{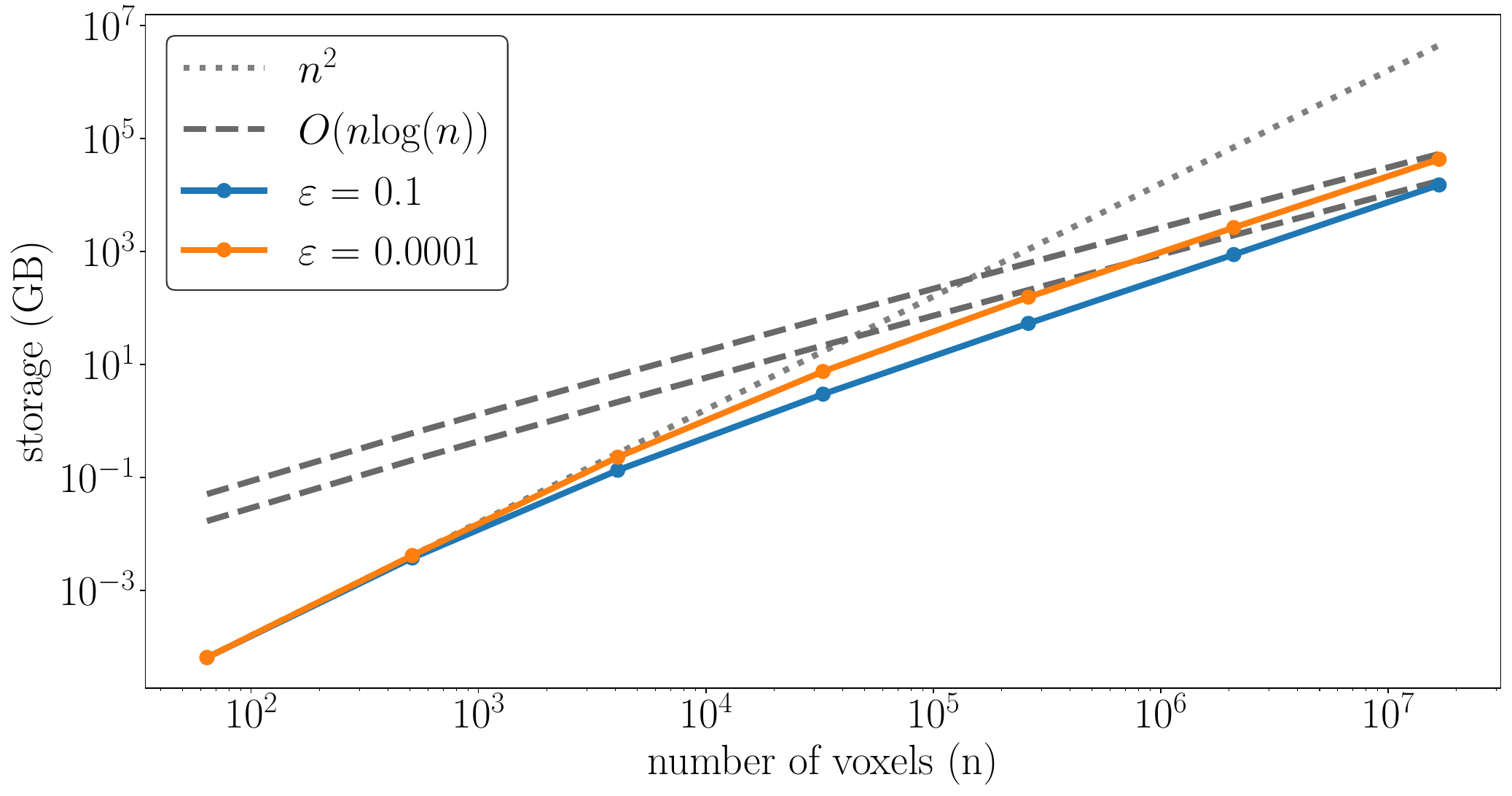}
    \caption{The memory requirements of the compression of matrix $B$ with tolerance $\varepsilon$ while using a fixed number of voxels per wavelength at multiple grids. The dashed grey lines depict an $\mathcal{O}(n\operatorname{log}(n))$ complexity with different constants.}
    \label{fig:complexity_VPWconst}
\end{figure}

\begin{table}[htbp]
\centering
\caption{The memory requirements of the compression of matrix $B$ with tolerance $\varepsilon$ while using a fixed number of voxels per wavelength at multiple grids. The cubic domain has sides measuring 1~cm. The memory for the dense matrix can be read in Table~\ref{tab:complexity_Fconst}.}
\label{tab:complexity_VPWconst}
\begin{tabular}{rrr|rr|rr}
\multicolumn{1}{l}{\multirow{2}{*}{\begin{tabular}[c]{@{}l@{}}number of\\ voxels ($n$)\end{tabular}}} &
  \multicolumn{1}{l}{\multirow{2}{*}{\begin{tabular}[c]{@{}l@{}}frequency\\ (kHz)\end{tabular}}} &
  \multicolumn{1}{l|}{\multirow{2}{*}{\begin{tabular}[c]{@{}l@{}}wavelength\\ (mm)\end{tabular}}} &
  \multicolumn{2}{c|}{$\varepsilon = 0.1$} &
  \multicolumn{2}{c}{$\varepsilon = 0.0001$} \\
\multicolumn{1}{l}{} &
  \multicolumn{1}{l}{} &
  \multicolumn{1}{l|}{} &
  \multicolumn{1}{l}{storage (GB)} &
  \multicolumn{1}{l|}{rate (\%)} &
  \multicolumn{1}{l}{storage (GB)} &
  \multicolumn{1}{l}{rate (\%)} \\ \hline
64         & 176.41    & 15.0     & 0.0000655  & 100.0 & 0.0000655  & 100.0 \\
512        & 352.82    & 7.5      & 0.00379  & 90.42 & 0.00419  & 100.0 \\
4,096      & 705.64    & 3.75     & 0.14      & 50.42 & 0.23      & 85.85 \\
32,768     & 1,411.27  & 1.875    & 3.0       & 17.49 & 7.58      & 44.13 \\
262,144    & 2,822.55  & 0.9375   & 53.80     & 4.89  & 156.76    & 14.26 \\
2,097,152  & 5,645.09  & 0.46875  & 887.58    & 1.26  & 2,665.72  & 3.79  \\
16,777,216 & 11,290.19 & 0.234375 & 15,164.12 & 0.34  & 43,089.31 & 0.96 
\end{tabular}
\end{table}

In contrast to the constant-frequency scenario, the results shown in Figure~\ref{fig:complexity_VPWconst} and Table~\ref{tab:complexity_VPWconst} are for an experiment in which the number of elements per wavelength is fixed as the frequency increases. The compressed matrix requires orders of magnitude less storage than the dense matrix. For the largest benchmark, the compression rate falls below 1\%, confirming the effectiveness of ACAPP+~DCC. However, the results did not reach an $\mathcal{O}(n \log (n))$ storage complexity in this experiment. This is consistent with observed behavior in related studies. For example, the analytical upper bounds on compression rank for frequency-dependent approximations are not as tight as for fixed frequencies or non-oscillatory kernels~\cite{engquist2018approximate}. Furthermore, other convergence studies with open-source BEM for P1-BIOs also did not achieve a log-linear complexity~\cite{betcke2017computationally}. Even though the $\mathcal{O}(n \log (n))$ storage complexity is not achieved in this experimental scenario, we emphasize that the compression rates are sufficiently good to run large-scale simulations of engineering interest that are impossible to perform without matrix compression.

\subsection{Transcranial ultrasound}

To evaluate the feasibility of the proposed compression algorithm in a realistic biomedical context, we applied the four described ACA variants to a VSIE model for transcranial ultrasound~\cite{ultrasonics2026}. The geometry and heterogeneous material parameters were derived from CT data\footnote{CT data obtained from: \href{https://www.morphosource.org/concern/media/000367572}{https://www.morphosource.org/concern/media/000367572}}. The grid consists of a skull slab containing a total of 199,693 rectangular cuboid voxels with dimensions of (0.489, 0.489, 0.5)~mm in Cartesian coordinates ($x$, $y$, $z$). We use a bowl transducer at a frequency of 500~kHz, with an exterior density of 1000~kg/m$^3$ and an exterior wave speed of 1500~m/s, reflecting the acoustic properties of water. The discretized VSIE matrix requires 972~GB of RAM for complex double-precision floating-point numbers. More information about CT data processing and acoustic configuration can be found in~\cite{ultrasonics2026}.

We used different octree partitionings for the boundary and volume operators. The surface octree had a minimum box size of 20 and the volumetric octree 30. The maximum octree depth was 4 in both instances. For the sustained convergence check, an $\alpha$ value of 2 was used.

\begin{figure}[htbp]
    \centering
    \includegraphics[width=.6\linewidth]{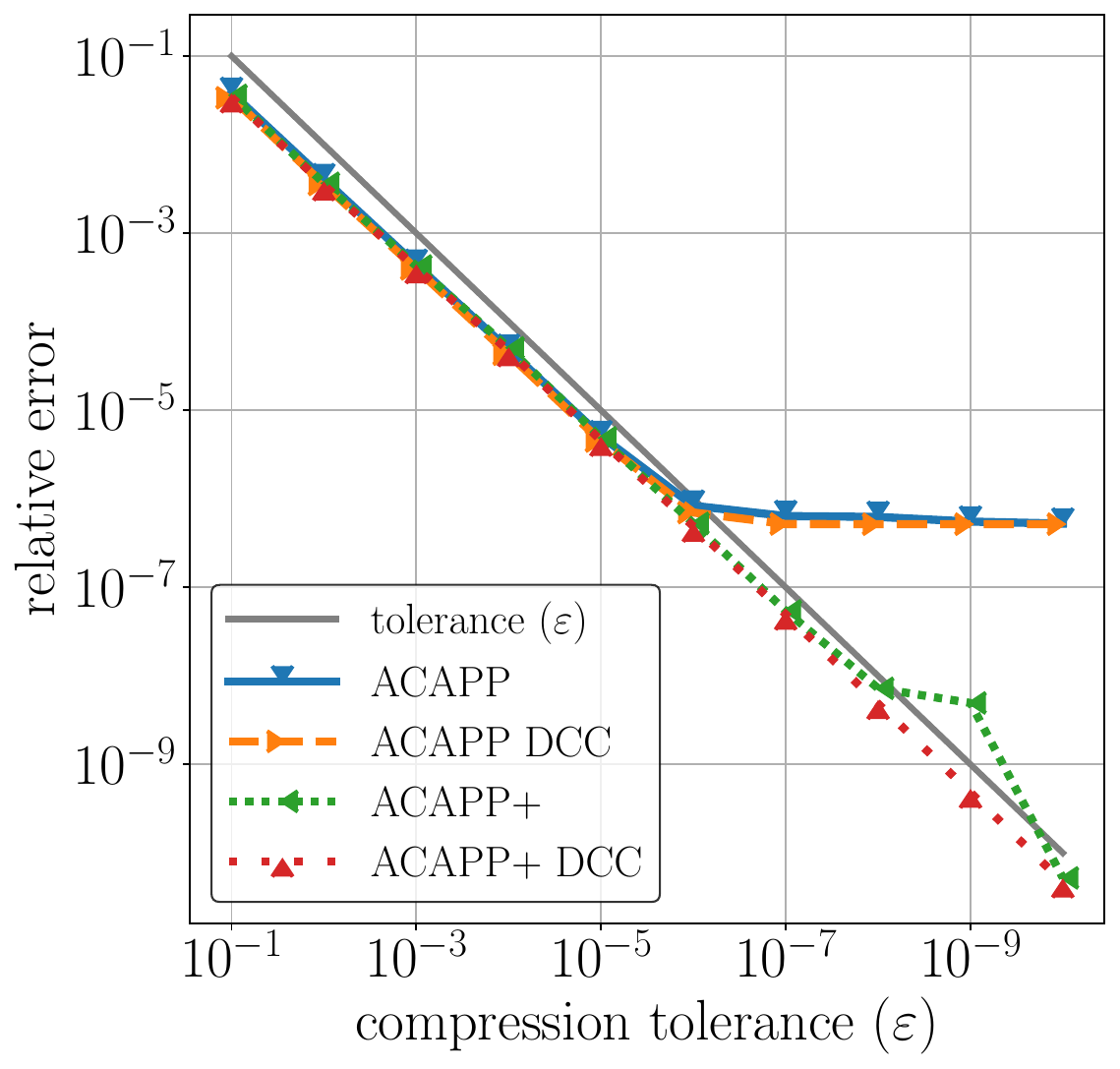}
    \caption{Relative errors obtained for the compression of the skull slab system with heterogeneous material parameters, and compression tolerances $\varepsilon\in\{10^{-1},10^{-2},\dots,10^{-10}\}$. The gray line represents the specified error tolerance for each compression.}
    \label{fig:skull_relative_errors}
\end{figure}

\begin{figure}[htbp]
    \centering
    \includegraphics[width=\linewidth]{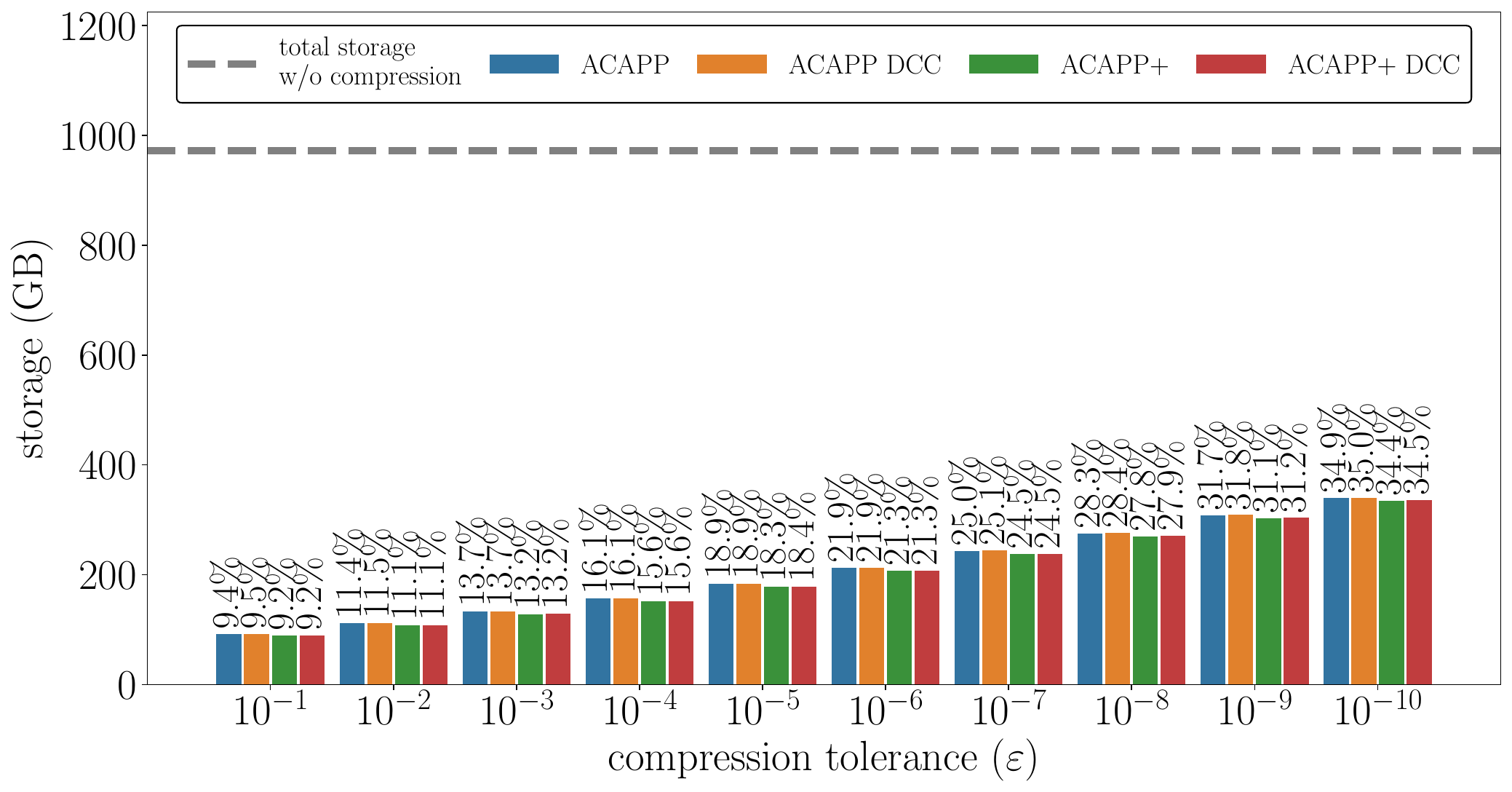}
    \caption{Memory storage and compression rates of the ACA algorithms, for the skull slab system with heterogeneous material parameters. The gray line depicts the storage of the original dense matrix.}
    \label{fig:skull_final_storages}
\end{figure}

Figure~\ref{fig:skull_relative_errors} shows the relative compression error~\eqref{eq:error_matvec} for the compressed VSIE on the skull slab. The results are consistent with the previous benchmarks. In particular, the double-layer operator has zero blocks since the skull slab has straight edges in the lateral directions. This causes early-convergence issues for various algorithms, as evidenced by stagnating accuracy for the ACAPP and ACAPP~DCC from $\varepsilon = 10^{-7}$ onwards. Additionally, the ACAPP+ fails to achieve the specified error tolerance at $\varepsilon = 10^{-9}$. Importantly, the ACAPP+~DCC algorithm always achieves the prescribed error margins.

Regarding the compression rates presented in Figure~\ref{fig:skull_final_storages}, the four algorithms again show similar performance. Even for the strict $\varepsilon = 10^{-10}$ tolerance, ACAPP+ DCC reduces the matrix storage by 65.5\% from 972~GB to approximately 335.34~GB.

\begin{figure}[htbp]
    \centering
    \includegraphics[width=\linewidth]{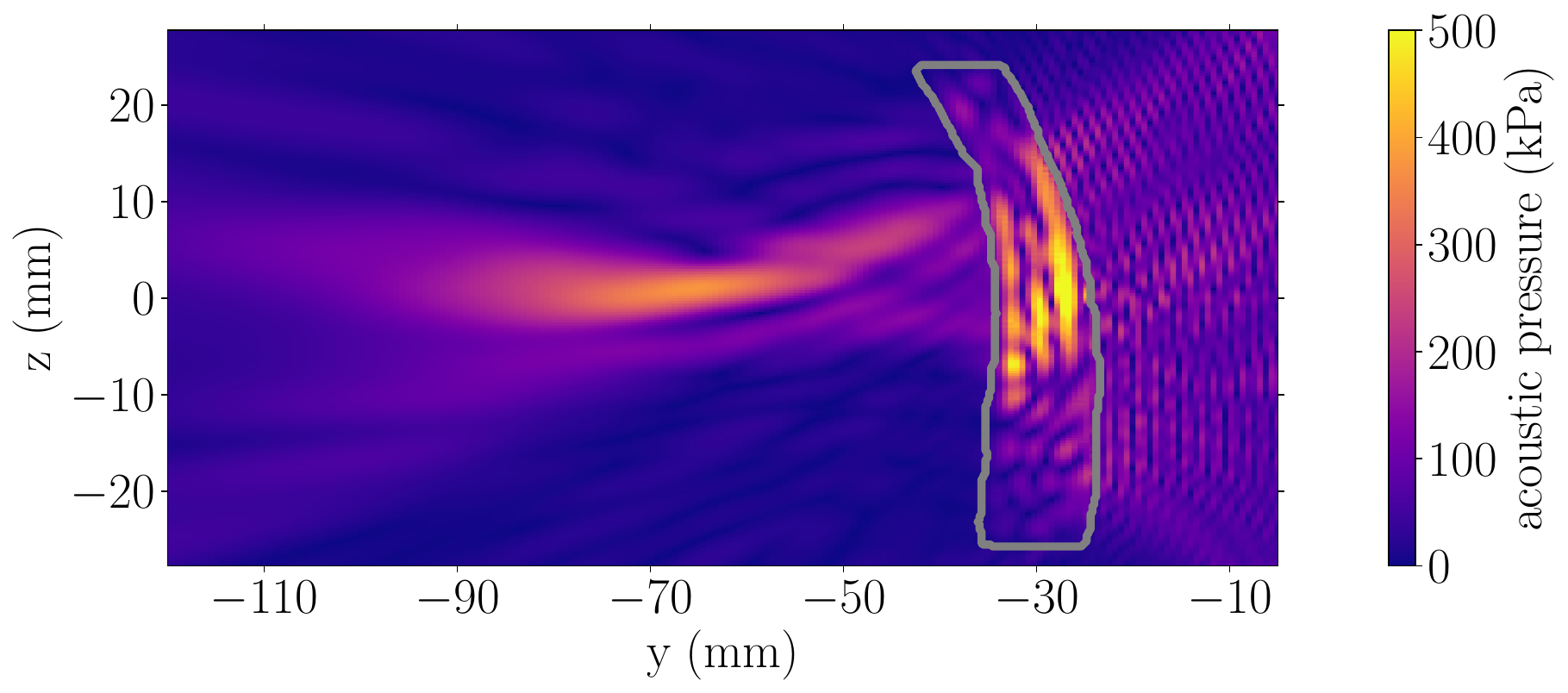}
    \caption{The magnitude of the acoustic field for a transcranial ultrasound simulation, on a slice at $x=0$. A bowl transducer at the global origin emits an acoustic field of 500~kHz in the negative $y$-direction. The exterior domain has a density of 1000~kg/m$^3$ and a wave speed of 1500~m/s, resembling water. The shape and material parameters of the skull bone are derived from CT data and stored in a voxel grid. The VSIE system was compressed with a tolerance of $\varepsilon = 10^{-5}$. The VSIE solution shows realistic wave patterns, with reflection on the right, energy absorption in the skull bone, and a transmitted focal region on the left. The gray line represents the contour of the skull bone.}
    \label{fig:skull_total_field}
\end{figure}

At a more common $\varepsilon = 10^{-5}$ tolerance, only 18.4\% of the memory is needed with ACAPP+ DCC compression. The solution of the compressed VSIE is visualized in Figure~\ref{fig:skull_total_field}. The pressure field inside the skull was solved using the VSIE, while the exterior pressure was calculated using the representation formula. This realistic acoustic field confirms the practical applicability of our robust methodology for large-scale acoustic simulations.

\section{Conclusions and future work}

Matrix compression is essential for solving integral equation formulations for realistic models, but early convergence remains a persistent challenge in both BIO and VIO frameworks. We proposed four algorithmic extensions to the standard ACAPP to achieve robust calculations. In our numerical simulations, the ACAPP+~DCC implementation is the only ACA version that successfully compresses all benchmarks within the prescribed error tolerances, thus avoiding the early-convergence problems of the alternative algorithms. Our new pivoting strategy was designed to maintain the log-linear storage complexity. In fact, computational experiments show that the ACAPP+~DCC variant achieves the same compression rates as traditional methods. Furthermore, the practical applicability of our methodology was confirmed by a large-scale transcranial ultrasound simulation with the VSIE. Hence, we conclude that our ACAPP+~DCC method provides a robust algorithm for performing realistic acoustic calculations using volume and boundary integral equations.

The presented robust compression algorithm successfully addressed early-convergence issues and reduced the memory requirements of acoustic integral equations. As future research, improved compression algorithms can be explored to further reduce memory consumption at large-scale high-frequency simulations. Efficient compression algorithms such as matrix recompression~\cite{borm2020hybrid}, $\mathcal{H}^2$-matrices~\cite{borm2010efficient}, and directional compression~\cite{borm2017approximation, borm2024memory} could be developed. They are expected to improve computational complexity, but their sensitivity to early-convergence problems must also be assessed.

\section*{Acknowledgments}

\subsection*{Software and data availability}
The Python code to reproduce the numerical results is available on GitHub (\href{https://github.com/aalmuna/RobustVSIEMatrixCompression}{github.com/aalmuna/RobustVSIEMatrixCompression}).

\subsection*{Declaration on generative AI}
During the preparation of this work, the authors used Google Gemini to improve the manuscript's language and readability. The authors reviewed and edited the content as needed and take full responsibility for the article's content.

\subsection*{Funding sources}
This work was financially supported by the Agencia Nacional de Investigación y Desarrollo (ANID), Chile, through FONDECYT 1230642 and DOCTORADO BECAS CHILE 72250237.

\printbibliography

\end{document}